\documentclass[11pt,twoside,reqno]{article}
\usepackage{setspace}
\usepackage{amssymb}
\usepackage{graphicx}
\usepackage{epstopdf}
\usepackage{color}

\usepackage{amsmath}

\usepackage{amsfonts}

\usepackage{amssymb}
\usepackage{graphics}
 \usepackage{bbm}
\usepackage{nomencl}

\usepackage{titlesec}
\usepackage[titletoc,toc,title]{appendix}

\titleformat{\section}{\large\bfseries}{\thesection.}{0.3em}{}
\titleformat{\subsection}{\bfseries}{\thesubsection.}{0.3em}{}
\titleformat{\subsubsection}{\small\bfseries}{\thesubsubsection.}{0.3em}{}

\makeglossary
 \allowdisplaybreaks[2]

\usepackage{latexsym}

\renewcommand{\skew}{\mathop{\rm skew}}
 \DeclareMathOperator{\sym}{sym}
\DeclareMathOperator{\tr}{tr} 
 \DeclareMathOperator{\dev}{dev}
\DeclareMathOperator{\Curl}{Curl\hskip.04truecm}
\DeclareMathOperator{\Div}{Div} 
\DeclareMathOperator{\inc}{inc}
\DeclareMathOperator{\sL}{\mathfrak{sl}}
\DeclareMathOperator{\so}{\mathfrak{so}}

\newcommand{\yieldlimit}{{\sigma}_{\mathrm{y}}}
\newcommand{\yieldspin}{{\widehat{\sigma}}_{\mathrm{y}}}

\newcommand{\yieldzero}{{\sigma}_{0}}
\newcommand{\yieldspino}{{\widehat{\sigma}}_{{0}}}

\newcommand{\C}{\mathbb{C}}

\newcommand{\PL}{p}

\newcommand{\BBR}{\mbox{$\mathbb{R}$}}

\newcommand{\SFH}{\mbox{$\mathsf{H}$}}
\newcommand{\SFL}{\mbox{$\mathsf{L}$}}
\newcommand{\SFQ}{\mbox{$\mathsf{Q}$}}
\newcommand{\SFV}{\mbox{$\mathsf{V}$}}
\newcommand{\SFW}{\mbox{$\mathsf{W}$}}

\newcommand{\SFZ}{\mbox{$\mathsf{Z}$}}

\newcommand{\bvarepsilon}{\mbox{$\bf\varepsilon$}}

\newcommand{\dsize}{\displaystyle}
\renewcommand{\div}{\mathop{\rm div}\nolimits}
\newcommand{\ba}{\mbox{\boldmath{$a$}}}
\numberwithin{equation}{section}

\newcommand{\norm}[1]{\lVert{#1}\rVert} 
\newcommand{\bfig}[2]{\begin{figure}\begin{center}\begin{picture}(341.8,#2)(
#1,0)}
\newcommand{\efig}[2]{\end{picture}\caption{#2.}\lbl{#1}\end{center}
\end{figure}}

\newcommand{\la}{\langle}
\newcommand{\ra}{\rangle}
\newcommand{\id}{{\boldsymbol{\mathbbm{1}}}}
\newcommand{\qed}{\qquad$\blacksquare$}
\newtheorem{theorem}{Theorem}[section]

\newtheorem{remark}{Remark}[section]

\newcommand{\be}{\begin{equation}}

\newcommand{\ee}{\end{equation}}


 \makeatletter
 \let\@fnsymbol\@arabic
 \makeatother
\begin{document}
\vskip-3truecm 
\title{
\vspace{-1.in} {\Large A canonical rate-independent model of
geometrically linear isotropic gradient plasticity with isotropic
hardening and plastic spin accounting for the Burgers
vector\texttt{}}}

\author{{\large Fran\c{c}ois Ebobisse{\footnote{Corresponding author, Fran\c{c}ois Ebobisse,
Department of Mathematics and Applied Mathematics, University of
Cape Town, Rondebosch 7700, South Africa, e-mail:
francois.ebobissebille@uct.ac.za}}\quad and\quad  Klaus
Hackl{\footnote {Klaus Hackl, Lehrstuhl f\"ur
Mechanik-Materialtheorie, Ruhr-Universit\"at Bochum,
Universit\"atsstrasse 150, 44801 Bochum, Germany, e-mail:
klaus.hackl@rub.de}} \quad and\quad Patrizio Neff{\footnote
{Patrizio Neff, Lehrstuhl f\"ur  Nichtlineare Analysis und
Modellierung, Fakult\"{a}t f\"ur Mathematik, Universit\"at
Duisburg-Essen, Thea-Leymann Str. 9, 45127 Essen, Germany, e-mail:
patrizio.neff@uni-due.de, http://www.uni-due.de/mathematik/ag\b{
}neff}} }\vspace{1mm}}



\date{\today}

 \maketitle

\begin{center}

 {{\small
\begin{abstract}
In this paper we propose a canonical variational framework for rate-independent phenomenological
geometrically linear gradient plasticity with plastic spin. The
model combines the additive decomposition of the total distortion
into non-symmetric elastic and plastic distortions, with a defect
energy contribution taking account of the Burgers vector through a dependence only on the dislocation density  tensor $\Curl p$ giving
rise to a non-symmetric nonlocal backstress, and isotropic hardening response only
depending on the accumulated equivalent plastic strain. The model is
fully isotropic and satisfies linearized gauge-invariance
conditions, i.e., only true state-variables appear. The model
satisfies also the principle of maximum dissipation which allows to show
existence for the weak formulation. For this result, a recently
introduced Korn's inequality for incompatible tensor fields is
necessary. Uniqueness is shown in the class of strong solutions.
For vanishing energetic length scale, the model reduces to classical
elasto-plasticity with symmetric plastic strain $\bvarepsilon_p$ and standard isotropic hardening.

   \end{abstract}}}
\end{center}
\noindent {\bf Key words:} plasticity, gradient plasticity, variational modeling, dissipation function, 
 geometrically necessary dislocations, incompatible distortions, rate-independent models,
isotropic hardening, generalized standard material, variational
inequality, convex analysis, associated flow rule, defect ener\-gy, dislocation density, plastic rotation, global dissipation inequality,
Burgers vector, plastic spin.
\vskip.2truecm\noindent {\bf AMS 2010 subject classification:}
35D30, 35D35, 74C05, 74C15, 74D10, 35J25.
\newpage
 {\small \tableofcontents}
\section{Introduction}\label{Intro}
Since the celebrated work of Tresca \cite{TRIESCA1872}, classical
plasticity has been cast within the years into a beautiful framework
in which both theoretical and computational aspects were examined
(see e.g. \cite{Lubliner2008, Martin1975, Alber, SIMOHUGHES,
Han-ReddyBook, Suquet1981, ddm}). Even perfect classical plasticity
has been recently revived by \cite{ddm, fragia2012,fragiamar2015}
with the
 use of the energetic approach for rate-independent processes developed by  \cite{Mielke2002, Mielke2005}.

On the other hand, a number of experimental results have shown
size-dependencies for the material behaviour in small scales
(micron/meso)  (see e.g. \cite{Fleckhutch1997,Stolkevans1998}).
However, classical plasticity models are scale independent and
therefore cannot capture those size-effects.
 This has led in the last thirty years to
an abundant literature (\cite{AIF1984, AIF1987, MUHAIF1991,
Fleckhutch2001, AIF2003, GURT2004, GUD2004, GURTAN2005, GURTAN2009,
FleckWillisI2009, FleckWillisII2009, REDDY2011}) on theories of
gradient plasticity with the aim of accommodating the experimentally
observed size effects mentioned above. The so-called {\it energetic
and dissipative length scales} have been involved. Moreover, effort
has also been made in the past years to provide mathematical results
for the initial boundary values problems and inequalities describing
some models of gradient plasticity (see for instance, \cite{DEMR1,
REM, EMR2008, NCA, EBONEFF, giacoluss,NESNEFF2012, NESNEFF2013,
ENR2015}). Several contributions on the computational aspects have
been made as well (\cite{DEMR2,NSW2009, BAREKLU2014, RWW2014}).

In most of the above-mentioned models of gradient plasticity, the
plastic rotation has been ignored. If a polycrystal is treated as a randomly oriented collection of grains,
it is clear that the plastic distortion $p$,
 which must then be seen as the average slip over all glide planes, will be non-symmetric.
  Therefore, plastic spin is a reality also in polycrystalline modelling.
  The situation is less clear when one aims at an overall effective phenomenological description
  in which individual glide planes are not resolved. It is possible to show that in a purely local
 isotropic theory the plastic spin can be suppressed without loss of generality. The situation is again different
in gradient-plasticity extensions, in which it is generally agreed that plastic spin is automatically
included (e.g. \cite{GURT2004}). However, no agreement has been reached on how to precisely include the effect of plastic spin.
 Our contribution aims at proposing a canonical framework to do exactly this.
 In \cite{GURT2004, Bardella2009,
Bardella2010, POHPEERLINGS2016} models discussing the role of the
plastic rotation have been proposed. For instance,
\cite{POHPEERLINGS2016} discusses the need to incorporate the
plastic rotation in an isotropic gradient plasticity framework in
order to capture some effects of a crystallographic model for a
large collection of grains in a polycrystal. In the mathematical
context, existence results for models with plastic spin have also
been obtained (\cite{NCA, EBONEFF, ENR2015}).

The modelling challenge which we faced in the past can be explained
as follows. Given the additive decomposition of the total
non-symmetric distortion (the displacement gradient $\nabla u$), is
it possible to write down a model with plastic spin (the plastic
distortion $p$ is not symmetric) and allow for a defect energy
depending on Nye's dislocation density tensor $\Curl p$ together with
an isotropic hardening response which is, however, only driven by
the accumulated equivalent plastic strain
$\gamma_p=\int_0^t\lVert\sym\dot{p}\rVert\,ds=\int_0^t\lVert\dot{\bvarepsilon}_p\rVert\,ds$,
and cast all that in the suitable convex variational framework of
the principle of maximum dissipation? In Section \ref{description}
we present exactly such a model. Our previous attempts of modelling
in this direction were based on the (rate-explicit) dual flow rule
but failed to satisfy the principle of maximum dissipation, \cite[p.
454]{GURTAN-BOOK}, see also \cite{hackl1997generalized,OrtRep99NEMD,CaHaMi02NCPM,hacklfischer}.\footnote{It is often assumed that
 the plastic evolution $\dot{p}$ associated
 with a state of yield  maximizes the dissipation relative to all admissible states. This is also equivalent to
 I'liushin's postulate (\cite{LUCCHSILHA1991}).}

The new model proposed in this paper, which involves only one
energetic length scale $L_c$ has some features which make it stand
out from other proposals in rate-independent gradient plasticity with plastic spin
as:
\begin{itemize}\item[$\bullet$] it allows for plastic spin in a most
transparent manner: for vanishing characteristic energetic length scale $L_c\to0$, the
plastic spin vanishes as well and the model turns into classical
elasto-plasticity with symmetric plastic strain $\bvarepsilon_p=\sym p$ and with isotropic hardening based only on the accumulated
equivalent plastic strain
$\gamma_p=\int_0^t\lVert\sym\dot{p}\rVert\,ds=\int_0^t\lVert\dot{\bvarepsilon}_p\rVert\,ds$;

\item[$\bullet$] it is completely isotropic and (linearized)
frame-indifferent;
\item[$\bullet$] it is (linearized) {\it gauge-invariant}: this means that it satisfies invariance
under compatible transformations of the reference system, i.e., in
the linearized context it is invariant under
$$\boxed{\begin{array}{llll}
 \nabla u(x) & \longrightarrow & \nabla u(x) +\nabla\vartheta_p(x)&\forall\vartheta_p\in C^2(\mathbb{R}^3,\mathbb{R}^3)\\
\,\,\,\,\,p(x)&\longrightarrow & \,\,\,\,\,\,p(x)+\nabla\vartheta_p(x) &\end{array}}\,,$$ which is also known as
translational T(3)-gauge invariance (\cite{LAZAR2000, LAZAR2002,
LAZANAS2008, ENS-new});

\item[$\bullet$] it contains only properly defined state-variables
(\cite{RUBIN2001, ENS-new}). In this context, notice that, as
mentioned in De Wit \cite[p.1478]{Dewit1981}: "\ldots the plastic
strain [$\sym p$] is not a state quantity, i.e., it cannot be
determined from the [current] state  of the body." Through a proper
definition of infinitesimal state-variables, this will be clearly
presented in \cite{ENS-new}.
\end{itemize}
\noindent

In this model, the hardening type response is depending on a
(nonlocal) kinematic term which is the non-symmetric backstress
contribution $\mu\,L^2_c\,\Curl \Curl p$, solely responsible for the
appearance of plastic spin or not and related to the geometrically
necessary dislocation (GND) density distribution. The
 isotropic  hardening is related to statistically stored
dislocations (SSD), which take into account a "plastically
homogeneous" effect as they accumulate already during a
macroscopically homogeneous deformation. Here, the SSD  evolution is
modelled by two isotropic hardening variables
$\gamma_p=\int_0^t\norm{\sym\dot{p}}\,ds$ and
$\omega_p=\int_0^t\norm{\skew\dot{p}}\,ds$. Hence, the full
plastic distortion, and not only its symmetric part, may contribute
to hardening. This is in accordance with the physical nature of
plastic flow since also the evolution of the skew-symmetric part of
$p$ indicates dislocation motion.  It is important to emphasize that no {\it spin cross-hardening} takes place in the proposed model, i.e., the situation where plastic flow in the
 plastic strain $\bvarepsilon_p=\sym p$ causes hardening in the plastic rotation
 evolution of $\skew p$ and vice-versa. This means that, in our model, only the
 accumulated equivalent plastic strain influences hardening in the
 evolution of the plastic strain and only the accumulated equivalent
 plastic rotation influences hardening in the evolution of the
 plastic rotation.

 It is noteworthy that classical linear Prager-type kinematical hardening cannot be accommodated in the "state-variable" approach adopted here since
  the corresponding backstress contribution $\bvarepsilon_p=\sym p$ as such is not a state-variable (see e.g. \cite{ONAT1996, RUBIN2001}).

Notwithstanding the use of the dislocation density tensor $\Curl p$,
we claim that our model is properly isotropic. In passing, notice
that taking $\Curl\bvarepsilon_p=\Curl\sym p$ is physically inadmissible since
 $\Curl\bvarepsilon_p$ is not a defect measure for $\bvarepsilon_p\in\mbox{Sym}(3)$. Rather, one should then
  take Kr\"oner's incompatibility tensor $\inc\bvarepsilon_p:=
 \Curl[(\Curl \bvarepsilon_p)^T]$. The possibilities to do exactly this will be explained   in the forthcoming paper \cite{ENS-new}.
 On the other hand, claims in the recent
literature \cite{STEIGUPTA}\footnote{Steigmann and Gupta \cite[p.410]{STEIGUPTA} put forward that:
"... the dislocation density [tensor] $\mbox{Curl}\, p$ is well-defined under symmetry transformations
only if the symmetry group is discrete." From that they conclude that it is not possible to obtain an isotropic
  plasticity model including $\Curl p$.}
 that dependencies of a model on the
dislocation density tensor $\Curl p$ exclude isotropy are also 
critically examined in \cite{ENS-new}.

 It is
sometimes argued that plastic spin is irrelevant in the case of
isotropy (\cite{KRISSTEIG2014}).\footnote{Krishnan and Steigmann \cite[p.722]{KRISSTEIG2014} argue that plastic spin
  associated with a flow rule for plastic evolution can be suppressed in the isotropic case without loss of generality. We understand that
   this is only true for the local theory, i.e., zero characteristic length $L_c=0$, as confirmed in \cite[p.511]{GURTAN-BOOK}.}
   The question whether one needs a theory with plastic spin is just the question whether one can work with
   a symmetric plastic strain tensor $\bvarepsilon_p$ as the only variable in a phenomenological plasticity theory. Our development clearly shows that
claims such as in \cite{{KRISSTEIG2014}} are unfounded and seem to indicate that there are different
notions involved of what isotropy precisely means. This subject is
 also  discussed further in \cite{ENS-new}.

 A remark concerning the mathematical treatment of single crystal
 plasticity is also in order. First, it is clear that such a theory is also a phenomenological model, albeit on a different
 scale.
  In the single crystal case the assumption of different glide
 systems lead to an immediate anisotropy of plastic flow and plastic spin is automatically included. However,
 the dislocation density contribution, when looked at it in detail,
 leads to a full gradient control of the plastic slip on each
 glide-plane. Therefore, the nonsymmetric plastic distortion $p$, which
 is the combined plastic slip  on each glide plane, is automatically
 controlled in the standard Sobolev space $H^1(\Omega)$
 (\cite{REDDY2011-II,BAREKLU2014}). By contrast, our isotropic framework means to
 give up detailed control of the plastic distortion due to additional 
 invariance conditions that have to be respected. The effect is that
 there is not even an immediate $H(\Curl)$-control of the plastic
 distortion. Therefore, the mathematically more challenging model
 is, without any doubt, the isotropic dislocation-based model with
 plastic spin treated here.

Notice that there are some similarities between our new isotropic model and
the early one proposed by Gurtin \cite{GURT2004}. In fact, both models
share: a complete isotropic formulation, decoupled evolution
equations into symmetric and skew-symmetric rates (isotropic
hardening possibly coming from both), a dissipation depending also on
plastic spin, the same defect and elastic energies, only an
energetic length scale connected to the dislocation density tensor
and both reduce to classical plasticity when the energetic length scale is
zero. Now, there are also nontrivial differences between the two
models. In fact, the model in \cite{GURT2004} is visco-plastic,
includes local nonsymmetric kinematical backstress due to dissipative
viscoplastic hardening, it is not cast into a variational framework and
does not have existence results so far. Also the model in \cite{GURT2004}
 involves a novel microforce balance as well as boundary conditions on the moving elastic-plastic boundary\footnote{
Unlike \cite{GURT2004}, no novel microforce balance needs to be
introduced in our model. Also, in our theory, nonstandard
(tangential) boundary
 conditions for the plastic distortion $p$ are always defined at the external boundary of the material only and the question
  on how to define them at a moving elastic-plastic boundary never arises. Any specific prescription of such boundary conditions
   at the elastic-plastic boundary could be in contradiction with the uniqueness result which we obtain for strong solutions.}
 and a dissipation function depending also on the gradient of the plastic distortion rate (see also \cite{QUOC2011}).
    The type of dissipation function considered in our model
 leads to an elastic region with Tresca-like branches and hence, in the flow rule in rate-explicit dual form, we get a case distinction to determine on which
  part of the yield surface the evolution takes place. In this, there are therefore similarities to crystal plasticity in which
   each glide plane has its own evolution and stresses are projected to the glide planes (see e.g. \cite{GURT2002}).
   In our model the non-symmetric Eshelby-type stress $\Sigma_E$ driving the plastic evolution
    is projected on $\sL(3)\cap\mbox{Sym}(3)$ (symmetric and traceless tensors) for the plastic strain evolution and $\so(3)$ (skew-symmetric tensors)
     for the plastic spin evolution.

Notice that the modelling capabilities of the model in \cite{GURT2004} have  been so far made relevant by many authors such as Bardella and co-authors \cite{Bardella2009, Bardegiacomini2008, Bardepante2015,PANTEBARDE2016} and also Poh and co-authors \cite{POHPEERLINGS2016}. So far, one still needs to consider a number of tests or examples to see whether the proposed model of isotropic hardening improves the results obtained in those papers or exhibits new features.

Let us emphasize that, while we will present the complete and rigorous mathematical existence theory to our model, the main thrust
 in this work is not only of analytical nature. It rather consists also in presenting that modeling framework for
  plastic spin which we deem to be the most suited one.

This paper is now structured as follows. In Section \ref{Notations},
we present some notations and definitions. In Section \ref{description}, we introduce
 various aspects of the model, in particular, the flow rule in both primal and dual
 formulations with the key role played by the dissipation function. In Sections
  \ref{alpha2-pos} and \ref{alpha2-zero}, we study mathematical aspects (existence and uniqueness) of the model while in Section \ref{limit-Lc}, we recover
   the classical plasticity framework when the characteristic length scale is set to be zero ($L_c\to0$).
\section{Some notational agreements and
definitions}\label{Notations} Let $\Omega$ be a bounded domain
 in $\BBR^3$ with Lipschitz continuous boundary $\partial\Omega$, which is occupied by the elastoplastic
body in its undeformed configuration. Let $\Gamma$ be a 
subset of $\partial\Omega$ with non-vanishing $2$-dimensional
Hausdorff measure. A material point in $\Omega$ is denoted by $x$
and the time domain under consideration is the interval $[0,T]$.\\
 For every $a,\,b\in\BBR^3$, we let $\la a,b\ra_{\BBR^3}$ denote the scalar
 product on $\BBR^3$ with associated vector
norm $|a|^2_{\BBR^3} = \la a, a\ra_{\BBR^3}$. We denote by
$\BBR^{3\times 3}$ the set of real $3\times 3$ tensors. The standard
Euclidean scalar product on $\BBR^{3\times 3}$ is given by $\la
A,\,B\ra_{\BBR^{3\times 3}} = \mbox{tr}\,\bigl[AB^T\bigr]$, where
$B^T$ denotes the transpose tensor of $B$. Thus, the Frobenius
tensor norm is $\norm{A}^2 = \la A,\,A\ra_{\BBR^{3\times 3}}$.
 In the following we omit the subscripts $\BBR^3$ and $\BBR^{3\times 3}$. The identity tensor on $\BBR^{3\times 3}$ will be denoted by
  $\id$, so that $\mbox{tr}(A) = \la A, \id\ra$.  The set
$\so(3):=\{X\in\BBR^{3\times 3}\,|\,\,X^T=-X\}$ is the Lie-Algebra
of skew-symmetric tensors.
 We let
$\mbox{Sym\,}(3):=\{X\in\BBR^{3\times 3}\,|\,\,X^T=X\}$ denote the
 vector space of symmetric tensors and $\sL(3):=\{X\in\BBR^{3\times
3}\,|\,\,\mbox{tr\,}(X)=0\}$ be the Lie-Algebra of traceless
tensors. For every $X\in\BBR^{3\times 3}$, we set
$\sym(X)=\frac12\bigl(X+X^T\bigr)$,
$\skew\,(X)=\frac12\bigl(X-X^T\bigr)$ and
$\dev(X)=X-\frac13\mbox{tr}\,(X)\,\id\in\sL(3)\,$ for the symmetric
part, the skew-symmetric part and the deviatoric part of $X$,
respectively. Quantities which are constant in space will be denoted
with an overbar, e.g., $\overline{A}\in\so(3)$ for the function
$A:\mathbb{R}^3\to\so(3)$ which is constant with constant value
$\overline{A}$.

The body is assumed to undergo infinitesimal deformations. Its
behaviour is governed by a set of constitutive
relations. Below is a list of variables and parameters used
throughout the paper:\begin{itemize}
\item[$\bullet$] $u$  is the displacement of the macroscopic material
points;
\item[$\bullet$] $\PL$ is the infinitesimal plastic distortion variable which is a
non-symmetric second order tensor, incapable of sustaining
volumetric changes; that is, $\PL\in\sL(3)$. The tensor $\PL$\,
represents the average plastic slip; $p$ is not a state-variable, while the rate $\dot{p}$ is;

\item[$\bullet$] $e=\nabla u -\PL$ is  the infinitesimal elastic distortion which  is a
non-symmetric second order tensor and is a state-variable;

\item[$\bullet$] $\bvarepsilon_p=\sym\, \PL$ is the symmetric infinitesimal plastic strain
tensor, which is also trace free, {$\bvarepsilon_p\in\sL(3)$;} $\bvarepsilon_p$ is not a state-variable;
the rate $\dot{\bvarepsilon}_p=\sym\dot{p}$ is a state-variable;

\item $\skew p$ is called plastic rotation or plastic spin;

\item[$\bullet$] $\bvarepsilon_e=\sym\,(\nabla u -\PL)$ is the symmetric infinitesimal  elastic
strain tensor and is a state-variable;
\item[$\bullet$] $\sigma$   is the Cauchy stress tensor which is a symmetric
second order tensor and is a state-variable;

\item[$\bullet$] $\yieldzero$ and $\yieldspino$  are the initial
yield stresses for plastic strain and plastic spin, respectively and both are state-variables;

\item[$\bullet$] $\yieldlimit$ and $\yieldspin$ are the current yield stresses
for plastic strain and plastic spin, respectively and both are state-variables;

\item[$\bullet$] $f$ is the body force;

\item[$\bullet$] $\Curl \PL=-\Curl e=\alpha$ is the dislocation density
tensor satisfying the so-called Bianchi identities $\Div\alpha=0$ and is a state-variable;

\item[$\bullet$] $\gamma_p=\dsize\int_0^t\lVert\sym\dot{\PL}\rVert\,ds=\int_0^t\lVert\dot{\bvarepsilon}_p\rVert\,ds$ is
the accumulated equivalent plastic
strain and is a state-variable;
\item[$\bullet$]$\omega_p=\dsize\int_0^t\lVert\skew\dot{\PL}\rVert\,ds$ is the accumulated equivalent plastic rotation and is a state-variable;
\item[$\bullet$]
$\dsize\int_0^t\hskip-.2truecm\sqrt{\dot{\gamma}_p^2+\dot{\omega}_p^2}\,ds=\int_0^t\norm{\dot{p}}\,ds$
 represents  the accumulated equivalent plastic distortion which is a state-variable.
\end{itemize}
\vskip.2truecm\noindent For isotropic media, the fourth order
isotropic elasticity tensor $\C_{\mbox{\scriptsize
iso}}:\mbox{Sym}(3)\to\mbox{Sym}(3)$ is given by
\begin{equation}
\C_{\mbox{\scriptsize iso}}\sym X = 2\mu\,\dev\,\sym X+\kappa
\,\tr(X) \id =2\mu\,\sym X+\lambda\,\tr(X)\id\label{C}
\end{equation}
for any second-order tensor $X$, where $\mu$ and $\lambda$ are the
Lam{\'e} moduli satisfying
\begin{equation}\label{Lame-moduli}
\mu>0\quad\mbox{ and }\quad 3\lambda +2\mu>0\,,
\end{equation} and $\kappa>0$ is the bulk modulus.
These conditions suffice for pointwise positive definiteness of the
elasticity tensor in the sense that there exists a constant $m_0 >
0$ such that
\begin{equation}
\forall X\in\mathbb{R}^{3\times 3}\mbox{ :}\,\quad\la \sym X,\C_{\mbox{\scriptsize iso}}\sym X\ra \geq m_0\,
\lVert\sym X\rVert^2\,. \label{ellipticityC}
\end{equation}

\vskip.2truecm\noindent The space of square integrable functions is
$L^2(\Omega)$, while the Sobolev spaces used in this paper are:
\begin{eqnarray}\label{sobolev-spaces}
\nonumber \mbox{H}^1(\Omega)&=&\bigl\{u\in
\mbox{L}^2(\Omega)\,\,|\,\,\mbox{grad\,}u\in
\mbox{L}^2(\Omega)\bigr\}\,,\qquad\quad\mbox{ grad}\,=\nabla\,,\\
\nonumber
&&\norm{u}^2_{H^1(\Omega)}=\norm{u}^2_{L^2(\Omega)}+\norm{\mbox{grad\,}u}^2_{L^2(\Omega)}\,,\qquad\forall u\in\mbox{H}^1(\Omega)\,,\\
 \mbox{H}(\mbox{curl};\Omega)&=&\bigl\{v\in
\mbox{L}^2(\Omega)\,\,|\,\,\mbox{curl\,}v\in
\mbox{L}^2(\Omega)\bigr\}\,,\qquad\mbox{curl\,}=\nabla\times\,,\\
\nonumber &&\norm{v}^2_{\mbox{\scriptsize
H}(\mbox{curl};\Omega)}=\norm{v}^2_{L^2(\Omega)}+\norm{\mbox{curl\,}v}^2_{L^2(\Omega)}\,,\,\,\quad\forall
v\in\mbox{H}(\mbox{curl;\,}\Omega)\,.
\end{eqnarray}
For every $X\in C^1(\Omega,\,\BBR^{3\times 3})$ with rows
$X_1,\,X_2,\,X_3$, we use in this paper the definition of $\Curl X$
in \cite{NCA, SVEN}:
\begin{equation}\label{def-Curl}\Curl X =\left(\begin{array}{l}\mbox{curl\,}X_1\\
\mbox{curl\,}X_2\\
\mbox{curl\,}X_3\end{array}\right)\in\BBR^{3\times
3}\,,\end{equation} for which $\Curl\,\nabla v=0$ for every $v\in
C^2(\Omega,\,\BBR^3)$. Notice that the definition of $\Curl X$
above is such that $(\Curl X)^Ta=\mbox{curl\,}(X^Ta)$ for every
$a\in\BBR^3$ and this clearly corresponds to the transpose of the
Curl of a tensor as defined in
\cite{GURTAN2005, GURTAN-BOOK}.\\

The following function spaces and norms will also be used later.
\begin{eqnarray}\label{Curl-spaces}
\nonumber \mbox{H}(\mbox{Curl};\,\Omega,\,\BBR^{3\times
3})&=&\Bigl\{X\in \mbox{L}^2(\Omega,\,\BBR^{3\times
3})\,\,\bigl|\,\,\mbox{Curl\,}X\in
\mbox{L}^2(\Omega,\,\BBR^{3\times 3})\Bigr\}\,,\\
\norm{X}^2_{\mbox{\scriptsize H}(\mbox{\scriptsize
Curl};\Omega)}&=&\norm{X}^2_{L^2(\Omega)}+\norm{\mbox{Curl\,}X}^2_{L^2(\Omega)}\,,\quad\forall
X\in\mbox{H}(\mbox{Curl;\,}\Omega,\,\BBR^{3\times 3})\,,\\
\nonumber
\mbox{H}(\mbox{Curl};\,\Omega,\,\mathbb{E})&=&\Bigl\{X:\Omega\to\mathbb{E}\,\,\bigl|\,\,X\in
\mbox{H}(\mbox{Curl};\,\Omega,\,\BBR^{3\times 3})\Bigr\}\,,
\end{eqnarray}
for $\mathbb{E}:=\sL(3)$ or
$\mbox{Sym}\,(3)\cap\sL(3)$.\vskip.2truecm\noindent
 We also consider the space
\begin{equation}\label{spacep-bc}\mbox{H}_0(\mbox{Curl};\,\Omega,\,\Gamma,\BBR^{3\times 3})\end{equation} as
the completion in the norm in  (\ref{Curl-spaces}) of the space
$\bigl\{q\in C^\infty(\Omega,\,\Gamma,\,\BBR^{3\times
3})\,\,\bigl|\,\,\,(q\times\,n)|_\Gamma=0\bigr\}\,.$ Therefore, this
space generalizes the tangential Dirichlet boundary condition
$$(q\times\,n)|_\Gamma=0\,$$
to be satisfied by the plastic distortion $\PL$ or the plastic
strain $\bvarepsilon_p:=\sym \PL$. The space
$$\mbox{H}_0(\mbox{Curl};\,\Omega,\,\Gamma,\mathbb{E})$$ is defined
as
 in (\ref{Curl-spaces}). \vskip.2truecm\noindent
 The divergence operator Div on second order
tensor-valued functions is also defined row-wise as
\begin{equation}\label{def-div}\mbox{Div}\,X=\left(\begin{array}{l}\mbox{div\,}X_1\\
\mbox{div\,}X_2\\
\mbox{div\,}X_3\end{array}\right)\,.\end{equation}

\section{The description of the model}\label{description}

\subsection{The balance equation} The conventional macroscopic
force balance leads to the equation of equilibrium
\begin{equation}
\div \sigma + f = 0\label{equil}
\end{equation}
in which $\sigma$ is the infinitesimal symmetric Cauchy stress and
$f$ is the body force.
 \subsection{Constitutive equations.} The
constitutive equations are obtained from a free-energy imbalance
together with a flow law that characterizes plastic behaviour. Since
the model under study involves plastic spin by which we mean that
the plastic distortion $p$ is not symmetric, we consider directly an additive
decomposition of the displacement gradient $\nabla u$ into elastic
and plastic components $e$ and $\PL$, so that
\begin{equation}
\nabla u = e + \PL\,, \label{displ-grad}
\end{equation}
with the nonsymmetric plastic distortion $\PL$ incapable of
sustaining volumetric changes; that is,
\begin{equation}
\tr(\PL)=\tr(\sym \PL)=\tr(\bvarepsilon_p) = 0\,. \label{trEp}
\end{equation}
 Here, $\bvarepsilon_e=\sym e=\sym (\nabla u-\PL)$ is the
infinitesimal elastic strain and $\bvarepsilon_p=\sym \PL$ is the
plastic strain while $\sym \nabla u=(\nabla
u+\nabla u^T)/2$ is the total strain.\\

 \noindent We consider  a free
energy in the additively separated form
\begin{eqnarray}\label{free-eng}
\Psi(\nabla u,p,\Curl\PL,\gamma_p,\omega_p):
&=&\underbrace{\Psi^{\mbox{\scriptsize lin}}_e(\sym
e)}_{\mbox{\small elastic energy}}\,\,
+\,\,\,\underbrace{\Psi^{\mbox{\scriptsize lin}}_{\mbox{\scriptsize
curl}}(\Curl \PL)}_{\mbox{\small defect
energy (GND)}}\\
\nonumber
&&\qquad\quad+\hskip-1.truecm\,\,\underbrace{\,\,\Psi_{\mbox{\scriptsize
 iso}}(\gamma_p,\omega_p)}_{\begin{array}{c}
\mbox{\small hardening energy (SSD)}\end{array}}\,,
 \end{eqnarray} where
 \begin{equation}\label{free-eng-expr}\begin{array}{ll} \Psi^{\mbox{\scriptsize
lin}}_e(\sym e):=\dsize\frac12\,\la\sym e,\C_{\mbox{\scriptsize iso}} \sym
e\ra, & \Psi^{\mbox{\scriptsize lin}}_{\mbox{\scriptsize
curl}}(\Curl \PL):=\dsize\frac12\,\mu\, L_c^2\norm{\Curl \PL}^2\,,\\&
\\ \Psi_{\mbox{\scriptsize
 iso}}(\gamma_p,\omega_p):=\dsize\frac12\,\mu\,\alpha_1|\gamma_p|^2
 +\dsize\frac12\,\mu\,\alpha_2|\omega_p|^2\,.&\end{array}\end{equation}
Here,  $L_c\geq0$ is an energetic length scale which characterizes the contribution of the defect energy density to the system, $\alpha_1>0$ and $\alpha_2\geq0$
are  nondimensional isotropic hardening constants, 
$\gamma_p$ and $\omega_p$ are isotropic hardening variables. The
defect energy is conceptually related to geometrically necessary
dislocations (GND). It is formed by the long-ranging stress-fields of excess dislocations and may be recovered by appropriate inelastic deformation.
The isotropic hardening energy is related to
statistically stored dislocations (SSD).\footnote{It is an easy
matter to generalize the defect-energy contribution as well as the
elasticity relation to the complete anisotropic setting. However,
this does not add anything to enhance understanding of the paper and
hence we leave these easy generalizations aside.} It is formed by the local stress-fields of all dislocations and can only be recovered in thermodynamical processes such as annealing, recrystallization or chemical reactions.

 \subsubsection{The
derivation of the dissipation inequality}\label{dis-ineq-section}
The local free-energy imbalance states that
\begin{equation}
\dot{\Psi} - \la\sigma,\dot{e}\ra - \la\sigma,\dot{\PL}\ra  \leq 0\
. \label{2ndlaw}
\end{equation}
Now we expand the first term, substitute (\ref{free-eng}) and get
\begin{equation}\label{exp-2ndlaw}
\la\C_{\mbox{\scriptsize iso}}\sym
e-\sigma,\sym\dot{e}\ra-\la\sigma,\dot{\PL}\ra+\mu
L_c^2\la\Curl\PL,\Curl\dot{\PL}\ra+\mu\,\alpha_1\gamma_p\,\dot{\gamma}_p+\mu\,\alpha_2\,\omega_p\,\dot{\omega}_p\leq0\,
,
\end{equation}
which, using arguments from thermodynamics gives the elastic
relation
\begin{equation}
\sigma = \C_{\mbox{\scriptsize iso}}\sym e=2\mu\, \sym(\nabla
u-\PL)+\lambda\, \tr(\nabla u-\PL)\id  \label{elasticlaw}
\end{equation}
and the local reduced dissipation inequality
\begin{equation}\label{reduced-diss}-\la\sigma,\dot{\PL}\ra+ \mu
L_c^2\la\Curl \PL,\Curl
\dot{\PL}\ra+\mu\,\alpha_1\gamma_p\,\dot{\gamma}_p+\mu\,\alpha_2\,\omega_p\,\dot{\omega}_p\leq0.\end{equation}
Now we integrate (\ref{reduced-diss}) over $\Omega$ and  get
\begin{eqnarray}\nonumber0&\geq&\int_{\Omega}\Bigl[-\la\sigma,\dot{\PL}\ra+ \mu
L_c^2\la\Curl \PL,\Curl \dot{\PL}\ra+\mu\,\alpha_1\gamma_p\,\dot{\gamma}_p+\mu\,\alpha_2\,\omega_p\,\dot{\omega}_p\Bigr]dx\\
\nonumber&=&-\int_{\Omega}\Bigl[\la\sigma,\dot{\PL}\ra+\mu
L_c^2\la\Curl\Curl \PL,\dot{\PL}\ra+\mu\,\alpha_1\gamma_p\,\dot{\gamma}_p+\mu\,\alpha_2\,\omega_p\,\dot{\omega}_p\\
&& \qquad+\sum_{i=1}^3\mbox{div}\Bigl(\mu
L_c^2\,\dot{\PL}^i\times(\Curl
\PL)^i\Bigr)\Bigr]\,dx\,.\label{reduced-diss2}
\end{eqnarray}
Using the divergence theorem we obtain
\begin{eqnarray}\label{reduced-diss3}
&&\nonumber\hskip-2truecm\int_{\Omega}\left[\la-\sigma+\mu L_c^2
\Curl\,\Curl
\PL,\dot{\PL}\ra+\mu\,\alpha_1\gamma_p\,\dot{\gamma}_p+\mu\,\alpha_2\,\omega_p\,\dot{\omega}_p\right]\,dx\\
&&\hskip3truecm+\,\sum_{i=1}^3\int_{\partial\Omega}\mu
L_c^2\langle\dot{\PL}^i\times(\Curl \PL)^i,\,n\rangle dS\leq0\, .
\end{eqnarray}In order to obtain a dissipation
inequality in the spirit of classical plasticity, we assume that the
infinitesimal plastic distortion $\PL$ satisfies the so-called {\it
linearized insulation condition}\footnote{Notice that the therminology ``insulation condition'' has been used by Polizzotto \cite{Polizzotto2009} and also in references therein.}
\begin{equation}\sum_{i=1}^3\int_{\partial\Omega}\mu
L_c^2\,\langle\dot{\PL}^i\times(\Curl \PL)^i,\,n\rangle
dS=0\,.\label{lin-sul}\end{equation} Under (\ref{lin-sul}) and
splitting the rates orthogonally in the scalar product
$\la\cdot\,,\cdot\ra$\,,
\begin{equation}\label{split-dotp}
\dot{\PL}=\sym\dot{\PL}+\skew\dot{\PL}\,,\end{equation} we then
obtain a global version of the reduced dissipation inequality\footnote{Gurtin \cite[p.4]{GURT2004} refers to Menzel and Steinmann
\cite{MENZSTEIN} and writes: "... but [they] satisfy the dissipation
 inequality [only] globally."}
\begin{eqnarray}
&&\nonumber\hskip-2truecm\int_{\Omega} [ \la\sigma
+\Sigma^{\mbox{\scriptsize lin}}_{\mbox{\scriptsize curl}},\dot
{\PL}\ra +g_1\,\dot{\gamma}_p+g_2\,\dot{\omega}_p]\,dx\geq0 \,
,\\
&\Leftrightarrow&\int_{\Omega} [ \la\sigma
+\sym\Sigma^{\mbox{\scriptsize lin}}_{\mbox{\scriptsize
curl}},\sym\dot {\PL}\ra + \la\skew\Sigma^{\mbox{\scriptsize
lin}}_{\mbox{\scriptsize curl}},\skew\dot {\PL}\ra
+g_1\,\dot{\gamma}_p+g_2\,\dot{\omega}_p]\,dx\geq0\,,
\,\label{diss-ineq}\end{eqnarray}

where \begin{equation}\label{stresses1}\Sigma^{\mbox{\scriptsize
lin}}_{\mbox{\scriptsize curl}}:=-\mu L_c^2\,\Curl\Curl \PL\,,
\qquad\quad g_1:=-\mu\,\alpha_1\,\gamma_p\,,\qquad\quad
g_2:=-\mu\,\alpha_2\,\omega_p\, .\end{equation}
For further use we define the non-symmetric Eshelby-type stress tensor driving the plastic evolution
\begin{equation}\label{Eshelby-stress}
\Sigma_E:= \sigma+\Sigma^{\mbox{\scriptsize
lin}}_{\mbox{\scriptsize curl}}\,,\end{equation}
with the non-symmetry relating only to the nonlocal term $\Sigma^{\mbox{\scriptsize
lin}}_{\mbox{\scriptsize curl}}$. In terms of $\Sigma_E$ the global reduced dissipation inequality can be expressed as
\begin{equation}\label{diss-ineq-Eshelby}
\int_{\Omega} [ \la\dev\sym\Sigma_E,\sym\dot {\PL}\ra + \la\skew\Sigma_E,\skew\dot {\PL}\ra
+g_1\,\dot{\gamma}_p+g_2\,\dot{\omega}_p]\,dx\geq0\,.
\end{equation}

The split used in
(\ref{split-dotp}) is a constitutive choice in that it will suggest
a suitable format on how to satisfy the inequality (\ref{diss-ineq})
in all deformation processes. In our previously proposed models (see
\cite{EBONEFF}), this split has not been used.
\subsubsection{The boundary conditions on the plastic
distortion}\label{bc-distortion}
 The condition (\ref{lin-sul}) is satisfied if
we
 assume for instance that the boundary is a perfect conductor. This means that the tangential component of $\PL$ vanishes on
$\partial\Omega$.
 In the context of dislocation dynamics these conditions express the requirement
  that there is no flux of the Burgers vector across a hard boundary.
Gurtin \cite{GURT2004} and also Gurtin and Needleman \cite{GURTNEED2005} introduce the
following different types of boundary conditions for the plastic distortion\begin{align}
   ( \dot{\PL}\times\,n)|_{\Gamma_{\rm hard}}&=0 \quad \text{"micro-hard" (perfect conductor)} \notag \\
    \dot{\PL}|_{\Gamma_{\rm hard}}&=0 \quad \text{"hard-slip"}\quad\mbox{(in the context of crystal plasticity)}\\
   ( \Curl \dot{\PL}\times\,n)|_{\Gamma_{\rm hard}}&=0 \quad \text{"micro-free"}\, . \notag
\end{align}
We specify a sufficient condition for the micro-hard boundary
condition, namely \begin{equation}\label{bc-plastic}
      ( \PL\times\,n)|_{\Gamma_{\rm hard}}=0
\end{equation}
and assume for simplicity only $\Gamma_{\rm
hard}=\partial\Omega=\Gamma$. Note that this boundary condition
constrains the plastic slip in tangential direction only, which is
what we expect to happen at the physical boundary $\Gamma_{\rm
hard}$.
\subsection{The flow rule}\label{flow-law}
\subsubsection{The flow rule in its primal formulation}\label{flow-law-primal}
Let $D:\mathbb{R}^2\to\mathbb{R}$ be  the function defined by
\begin{equation}\label{pre-diss}
D(s,t):=\sqrt{\yieldzero^2\,s^2+\yieldspino^2\,t^2}\,,
\end{equation}
 where $\yieldzero,\,\yieldspino>0$ are the initial yield stresses for symmetric strain $\sym p$ and skew-symmetric spin $\skew p$, respectively.\footnote{Both values together will define the
 elastic domain in the stress space and this domain must have nonempty interior. Therefore, we need $\yieldzero,\,\yieldspino>0$. Without
  isotropic hardening the elastic domain turns out to be \\$\{\Sigma_E\in\mathbb{R}^{3\times 3}\,\,|\,\,\norm{\dev\sym\Sigma_E}\leq\yieldzero,\,\,
  \norm{\skew\Sigma_E}\leq\yieldspino
  \}$.}\vskip.1truecm\noindent
  We consider the
dissipation function $\Delta$ defined by\,\footnote{ Gurtin
\cite[p.2554]{GURT2004} notes: "One would expect that, plastically,
the material response to spin differs to straining, and that
straining and spin each incur dissipation." Gurtin's choice of the
dissipation function in \cite{GURT2004} corresponds to $\yieldspino=\chi\,\yieldzero\geq0$ in (\ref{pre-diss}).
Also, Gurtin  \cite[p.2558]{GURT2004} takes
$\chi\to0$ formally and recovers  classical plasticity. If we want to take $\yieldspino\to0$ in our setting, then we
 encounter a problem described in Section \ref{hatsigma-equalzero}.}
\begin{equation}\label{diss-funct1}
\Delta(q,\eta,\beta):=\left\{\begin{array}{ll}\hskip-.2truecm
D(\lVert\sym q\rVert,\lVert\skew q\rVert) & \mbox{if }\norm{\sym
q}\leq\eta\quad\mbox{and}\quad
\norm{\skew q}\leq\beta\,,\\\\
\hskip-.2truecm\infty & \mbox{otherwise}\,.\end{array}\right.
\end{equation} 
The flow rule in its primal
formulation can be derived using the principle of the minimum of the dissipation function \cite{hackl1997generalized,OrtRep99NEMD,CaHaMi02NCPM}, stating that the rate of the internal variables is the minimizer of a functional $L$ consisting of the sum of the rate of the free energy and the dissipation function with respect to appropriate boundary conditions,
\begin{equation}\label{min-diss}
L = \int_\Omega [ \dot\Psi + \Delta ] dx.
\end{equation}
The principle of the minimum of the dissipation function is closely related to the principle of maximum dissipation. Both are not physical principles but thermodynamically consistent selection rules which turn out to be convenient if no other information is available or if existing flow rules are to be extended to a more general situation. For a detailed investigation, see \cite{hacklfischer}. A very general exposition for coupled physical processes is worked out in \cite{hackl2011,hackl2011a}. Applications to the evolution of plastic microstructures can be found in \cite{hackl2012model,hackl2008relaxed,hackl2014variational,koch}.
\vskip.3truecm\noindent
Employing a partial integration, the stationarity conditions of (\ref{min-diss}) can be compactly stated as
\begin{equation}\label{flow-law-primal1}
\Sigma_p\in\partial\Delta(\dot{\Gamma}_p)\qquad\mbox{ where
}\Sigma_p=(\sigma+\Sigma_{\mbox{\scriptsize
curl}}^{\mbox{\scriptsize lin}},g_1,g_2)\,\quad\mbox{ and }\quad
\Gamma_p=(\PL,\gamma_p,\omega_p)\,,\end{equation} and where
$\partial\Delta$ denotes the subdifferential of $\Delta$. That is,
for $\Sigma_p\in\partial\Delta(\dot{\Gamma}_p)$ we must have
\begin{eqnarray}\label{flow-law-primal2}
\nonumber\Delta(\overline{\Gamma})&\geq &\Delta(\dot{\Gamma}_p)
+\la\Sigma_p,\overline{\Gamma}-\dot{\Gamma}_p\ra\\
\nonumber&=&\Delta(\dot{\Gamma}_p)+\la\sigma+\Sigma_{\mbox{\scriptsize
curl}}^{\mbox{\scriptsize lin}},q-\dot{\PL}\ra
+g_1\,(\eta-\dot{\gamma}_p)+g_2\,(\beta-\dot{\omega}_p)\\
&=&\Delta(\dot{\Gamma}_p)+\la\Sigma_E,q-\dot{\PL}\ra
+g_1\,(\eta-\dot{\gamma}_p)+g_2\,(\beta-\dot{\omega}_p)\,,
\end{eqnarray} for every $\overline{\Gamma}=(q,\eta,\beta)$. By choosing $\overline{\Gamma}=(0,0,0)$ in (\ref{flow-law-primal2}), we get
the reduced dissipation inequality in pointwise form
\begin{eqnarray}
&&\nonumber\hskip-2truecm \la\Sigma_E,\dot
{\PL}\ra +g_1\,\dot{\gamma}_p+g_2\,\dot{\omega}_p\,\geq0 \,
,\\
&\Leftrightarrow&\la\dev\sym\Sigma_E,\sym\dot {\PL}\ra + \la\skew\Sigma_E,\skew\dot {\PL}\ra
+g_1\,\dot{\gamma}_p+g_2\,\dot{\omega}_p]\,\geq0
\,\label{diss-ineq-pointwise}\end{eqnarray}

\subsubsection{The flow rule in its dual
formulation}\label{flow-law-dual} 
\begin{figure}[h!]\tiny
\centering

\includegraphics[width=0.6\textwidth]{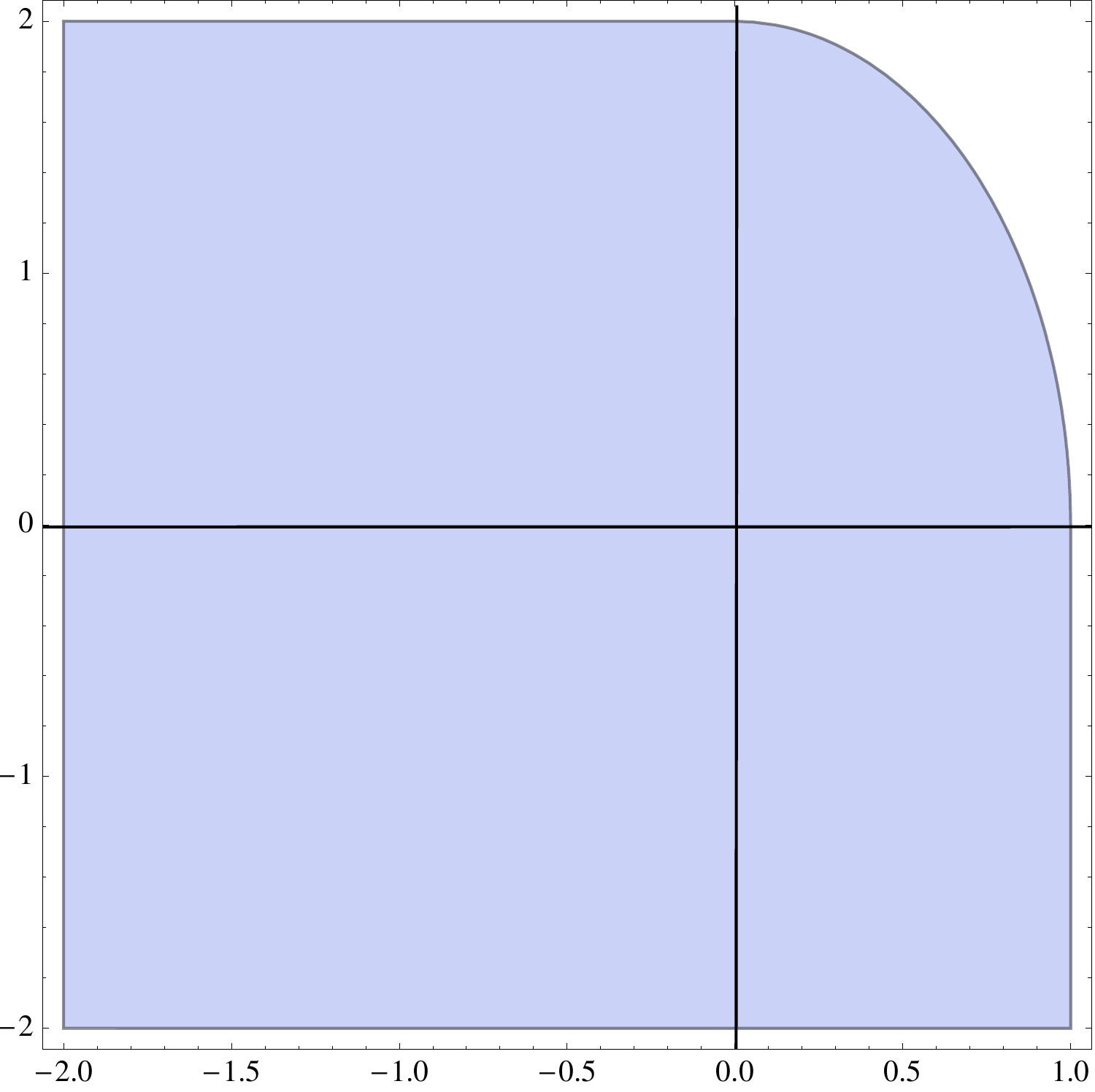}
\put(2,63){$\mathcal{S}_1$}
\put(-167,272){$\mathcal{S}_2$}\put(-6,135){$\sigma_0$}\put(-12,198){$\mathcal{S}_3$}\put(-88,266){$\hat{\sigma}_0$}
\put(-95,133){$0$} \caption{The set $\mathcal{K}$ in the $AB$-plane.
On $\mathcal{S}_1$ the flow is only driven by the symmetric rate part
$\sym\dot{p}$ as $\skew\dot{p}=0$. On $\mathcal{S}_2$ the flow is
only driven by the skew-symmetric rate part $\skew\dot{p}$ as
$\sym\dot{p}=0$. On $\mathcal{S}_3$ the flow is driven by both
symmetric and skew-symmetric rate parts $\skew\dot{p}$ and $\sym\dot{p}$.
}\label{fig:the-set-K}
\end{figure}
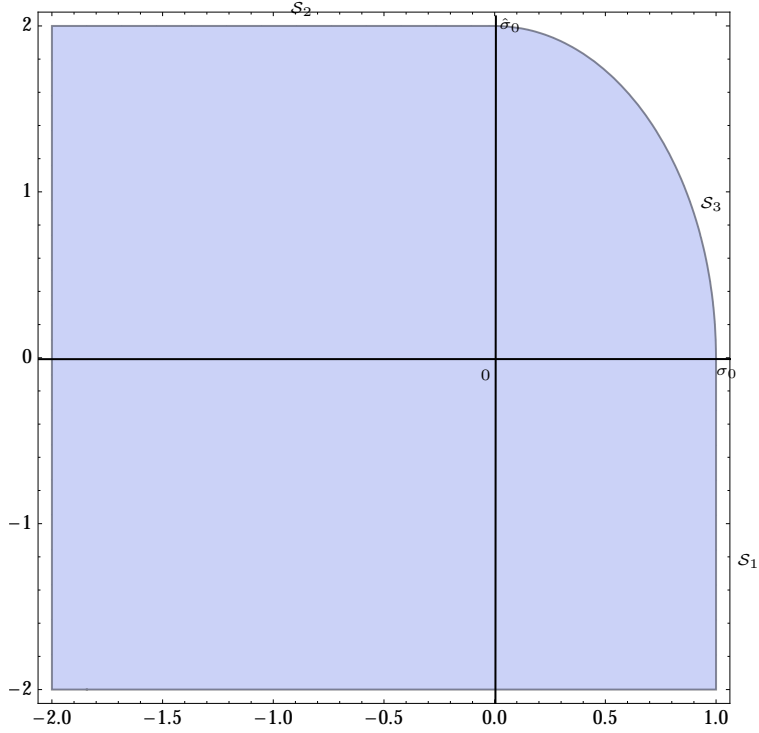

While the flow rule  in the primal
formulation is extremely condensed and will allow us a mathematical
treatment (existence), we need the representation of the flow rule
in the dual formulation in most computational implementations and
for the uniqueness proof in Section \ref{uniqueness}. For this
formulation of the flow rule we need to derive the set of admissible
(generalized) stresses $\mathcal{E}$ (the elastic domain)
corresponding to the dissipation function $\Delta$. According to the
principle of maximum dissipation,\footnote{which, again, is not a principle,
but a useful and often made simplifying assumption.} the flow rule in
dual form is formulated in the context of convex analysis as
\begin{equation}\label{flow-law-dual1}
\dot{\Gamma}_p\in N_{\mathcal{E}}(\Sigma_p)\quad
\Leftrightarrow\quad
\la\dot{\Gamma}_p,\overline{\Sigma}-\Sigma_p\ra\leq0\qquad\forall\,\overline{\Sigma}\in\mathcal{E}\,,\end{equation}
where $N_{\mathcal{E}}(\Sigma_p)$ is the normal cone to the set
$\mathcal{E}$ of admissible stresses at $\Sigma_p$. Therefore, we
need  to find the set $\mathcal{E}$. In the context of convex
analysis, the indicator function $I_\mathcal{E}$ of the set
$\mathcal{E}$ is the Fenchel-Legendre conjugate of the dissipation
function $\Delta$. Let us find the set $\mathcal{E}$ whose interior
$\mbox{int}(\mathcal{E})$ is the elastic domain and its boundary
$\partial\mathcal{E}$ is the yield surface. \vskip.2truecm\noindent
For $\Sigma_p=(\Sigma_E,g_1,g_2)$ with $\Sigma_E:=\sigma
+\Sigma_{\mbox{\scriptsize curl}}^{\mbox{\scriptsize lin}}$, we have
\begin{eqnarray}\label{indic-K1}
\nonumber I_\mathcal{E}(\Sigma_p)&=&\sup\bigl\{\la\Sigma_p,\Gamma\ra-\Delta(\Gamma)\,\,\,|\,\,\,\Gamma=(q,\eta,\beta)\bigr\}\\
\nonumber&=&\sup\bigl\{\la\Sigma_E,q\ra+g_1\,\eta+g_2\,\beta-\Delta(q,\eta,\beta)\,\,|\,\,\norm{\sym
q}\leq\eta\,,\,\,\,\norm{\skew q}\leq\beta\bigr\}\\
\nonumber&=&\sup_q\Bigl[\sup_{\eta,\,\beta}\bigl\{\la\Sigma_E,q\ra+g_1\,\eta+g_2\,\beta-\Delta(q,\eta,\beta)\,\,\,|\,\,\,\norm{\sym
q}\leq\eta\,,\,\,\,\norm{\skew q}\leq\beta\bigr\}\Bigr]\\
\nonumber&=&\sup_q\bigl\{\la\Sigma_E,q\ra+g_1\,\norm{\sym\,q}+g_2\,\norm{\skew\,q}-\Delta(q,\norm{\sym q},\norm{\skew q})\,\bigr\} \\
&=&\sup_q\left\{\hskip-.2truecm\begin{array}{l}\la\dev\sym\Sigma_E,\sym
q\ra+\la\skew\Sigma_E,\skew q\ra\\
\qquad+\,\,g_1\,\norm{\sym q}+g_2\,\norm{\skew
q}-\Delta(q,\norm{\sym q},\norm{\skew
q})\end{array}\hskip-.2truecm\right\}\,,
\end{eqnarray} where the supremum with respect to $\eta$ and $\beta$
 is achieved for $\eta=\norm{\sym q}$ and $\beta=\norm{\skew q}$ since
$g_1\leq 0$ and $g_2\leq0$.
\vskip.1truecm\noindent
Now taking the supremum with respect to $q$ and using the fact that
$\la\Sigma_E,q\ra$ is maximum with respect to $q$ only when $q$ is
in the direction of $\Sigma_E$, we find that it is not restrictive
to assume that
\begin{equation}\label{max-dotprod}\sym\,q=s\,\frac{\dev\,\sym\,\Sigma_E}{\norm{\dev\,\sym\,\Sigma_E}}\quad\mbox{
and }\quad
\skew\,q=t\,\frac{\skew\,\Sigma_E}{\norm{\skew\,\Sigma_E}}\,.\end{equation}
We then obtain
\begin{equation}\label{indic-K2} I_{\mathcal{E}}(\Sigma_p)=
\sup_{s\geq0,\,t\geq0}\Bigl\{ s\,(\norm{\dev\,\sym\Sigma_E}+g_1)+
t\,(\norm{\skew\Sigma_E}+g_2)-\sqrt{\yieldzero^2\,s^2+\yieldspino^2\,t^2}\,\Bigr\}\,.\end{equation}
To simplify the function of $s$ and $t$ to be maximized in
(\ref{indic-K2}), we set
\begin{equation}\label{new-param}
A:=\norm{\dev\,\sym\Sigma_E}+g_1\quad\mbox{ and }\quad
B:=\norm{\skew\Sigma_E}+g_2\,,
\end{equation}
 and hence,
\begin{equation}\label{indic-K3}I_{\mathcal{E}}(\Sigma_p)=
\sup_{s\geq0,\,t\geq0}\Bigl\{A\, s+
B\,t-\sqrt{\yieldzero^2\,s^2+\yieldspino^2\,t^2}\,\Bigr\}\,.\end{equation}
Notice immediately that 
\begin{eqnarray}\label{indic-K4}\nonumber I_{\mathcal{E}}(\Sigma_p) &=&
\sup_{s\geq0,\,t\geq0}\left\{\frac{A}{\yieldzero} s+
\frac{B}{\yieldspino}t-\sqrt{s^2+t^2}\,\right\}\,\\
&=&\left\{\begin{array}{ll}0
&\mbox{if }
\left\{\begin{array}{ll}A\leq \yieldzero &\mbox{ if }\quad \,B\leq 0\\&\\
B\leq \yieldspino &\mbox{ if }\quad \,A\leq0\\&\\
\dsize\frac{A^2}{\yieldzero^2}+\dsize\frac{B^2}{\yieldspino^2}\leq1
&\mbox{ if }\left\{\begin{array}{l}A\geq0\\B\geq0\end{array}\right.
\end{array}\right.\\&\\
\infty &\mbox{ otherwise}\,.\end{array}\right.
\end{eqnarray}  Let us
now introduce a set $\mathcal{K}\subset\mathbb{R}^2$ needed for
elucidating the branching behaviour of our flow rule and defined by
\begin{equation}\label{K}
\mathcal{K}:=\mathcal{K}_1\cup\mathcal{K}_2\cup\mathcal{K}_3\,,\end{equation}
where
\begin{eqnarray*}\label{K1}\mathcal{K}_1&=&\bigl(-\infty,\,\sigma_0\bigr]\times\bigl(-\infty,\,\,0\bigr]\,,
\qquad\mathcal{K}_2\,=\,\bigl(\,-\infty,\,\,0\,\bigr]\,\times(-\infty,\widehat{\sigma}_0]\,,\\
\label{K3}\mathcal{K}_3&=&\left\{(A,B)\in\mathbb{R}^+\times\mathbb{R}^+\,\,\Bigl|\,\,\,
\dsize\frac{A^2}{\yieldzero^2}+\dsize\frac{B^2}{\yieldspino^2}\leq1\right\}\,.\end{eqnarray*}
 The set $\mathcal{K}$ in the $AB$-plane is represented
graphically in Figure \ref{fig:the-set-K}. Notice that the set $\mathcal{K}$ itself is {\bf not} the elastic domain.
In our setting, the elastic domain is then defined as the interior of the set
\begin{equation}\label{the-setK3}
\mathcal{E}=\Bigl\{(\Sigma_E,g_1,g_2)\in\mathbb{R}^{3\times
3}\times\mathbb{R}^-\times\mathbb{R}^-\,\,|\,\,(\norm{\dev\,\sym\Sigma_E}+g_1,\norm{\skew\Sigma_E}+g_2)\in\mathcal{K}\Bigr\}\,.
\end{equation} In other terms, the set $\mathcal{E}$, which is also called the set of admissible stresses, is expressed as
\begin{equation}\label{the-setK4}
\mathcal{E}=\mathcal{E}_1\cup\mathcal{E}_2\cup
\mathcal{E}_3\,,\end{equation} where

\begin{eqnarray*}\label{E1}\mathcal{E}_1&=&\Bigl\{(\Sigma_E,g_1,g_2)\,\,\,|\,\,\,\norm{\dev\,\sym\Sigma_E}\leq
-g_1+\yieldzero\,,\quad
 \,\norm{\skew\Sigma_E}\leq -g_2\hskip.8truecm\Bigr\}\,,\\
\label{E2}\mathcal{E}_2&=&\Bigl\{(\Sigma_E,g_1,g_2)\,\,\,|\,\,\,\norm{\dev\,\sym\Sigma_E}\leq
-g_1\,,\qquad\,\,\quad\norm{\skew\Sigma_E}\leq -g_2+\yieldspino\Bigr\}\,,\\
\label{E3}\mathcal{E}_3&=&\left\{\hskip-.2truecm\begin{array}{l}(\Sigma_E,g_1,g_2)\,\,\,|\,\,\,\norm{\dev\,\sym\Sigma_E}\geq
-g_1,\qquad\quad\,\,\norm{\skew\Sigma_E}\geq-g_2,\\\\
\qquad\mbox{and }\hskip1.1truecm\dsize\frac{(\norm{\dev\,\sym\Sigma_E}+g_1)^2}{\yieldzero^2}+
\dsize\frac{(\norm{\skew\Sigma_E}+g_2)^2}{\yieldspino^2}\leq1\end{array}\hskip.3truecm\right\}\,.\end{eqnarray*}
 Hence, the yield surface is given by
\begin{equation}\label{Yield-surface}
\partial\mathcal{E}=\mathcal{S}_1\cup\mathcal{S}_2\cup\mathcal{S}_3\end{equation}
with

\begin{eqnarray}\nonumber\mathcal{S}_1&=&\Bigl\{(\Sigma_E,g_1,g_2)\,\,\,|\,\,\,\norm{\dev\,\sym\Sigma_E}=
-g_1+\yieldzero\,,\quad
 \,\norm{\skew\Sigma_E}\leq -g_2\hskip.8truecm\Bigr\}\,,\\
\label{S2}\mathcal{S}_2&=&\Bigl\{(\Sigma_E,g_1,g_2)\,\,\,|\,\,\,\norm{\dev\,\sym\Sigma_E}\leq
-g_1\,,\qquad\,\,\quad\norm{\skew\Sigma_E}= -g_2+\yieldspino\Bigr\}\,,\\
\nonumber\mathcal{S}_3&=&\left\{\hskip-.2truecm\begin{array}{l}(\Sigma_E,g_1,g_2)\,\,\,|\,\,\,\norm{\dev\,\sym\Sigma_E}\geq
-g_1,\qquad\quad\,\,\norm{\skew\Sigma_E}\geq-g_2,\\\\
\qquad\mbox{and }\hskip1.1truecm\dsize\frac{(\norm{\dev\,\sym\Sigma_E}+g_1)^2}{\yieldzero^2}+
\dsize\frac{(\norm{\skew\Sigma_E}+g_2)^2}{\yieldspino^2}=1\end{array}\hskip.3truecm\right\}\,.\end{eqnarray}

\begin{remark}\label{rem-dissip-01}{\rm 
We could consider a more general dissipation function corresponding e.g. to the function
\begin{equation}\label{pre-dissgen}
\widehat{D}(s,t):=r_1\,s+r_2\,t+\sqrt{\yieldzero^2\,s^2+\yieldspino^2\,t^2}\,\mbox{ \qquad
with $r_1,r_2\geq0$}\,.
\end{equation}
For such a choice, we get by easy calculations, the set of admissible stresses 
\begin{equation}\label{proof4}
\overline{\mathcal{E}}=\Bigl\{(\Sigma_E,g_1,g_2)\in\mathbb{R}^{3\times
3}\times\mathbb{R}^-\times\mathbb{R}^-\,\,|\,\,(\norm{\dev\,\sym\Sigma_E}+g_1+r_1,\norm{\skew\Sigma_E}+g_2+r_2)\in\mathcal{K}\Bigr\}\,,
\end{equation}
which corresponds to a dilation of the set in (\ref{the-setK3}).
Such a choice will not add any particular feature to the current model. In fact, 
 this simply corresponds to the expansion of
the initial elastic domain (i.e. before isotropic hardening takes place).}\end{remark}
\vskip.3truecm\noindent

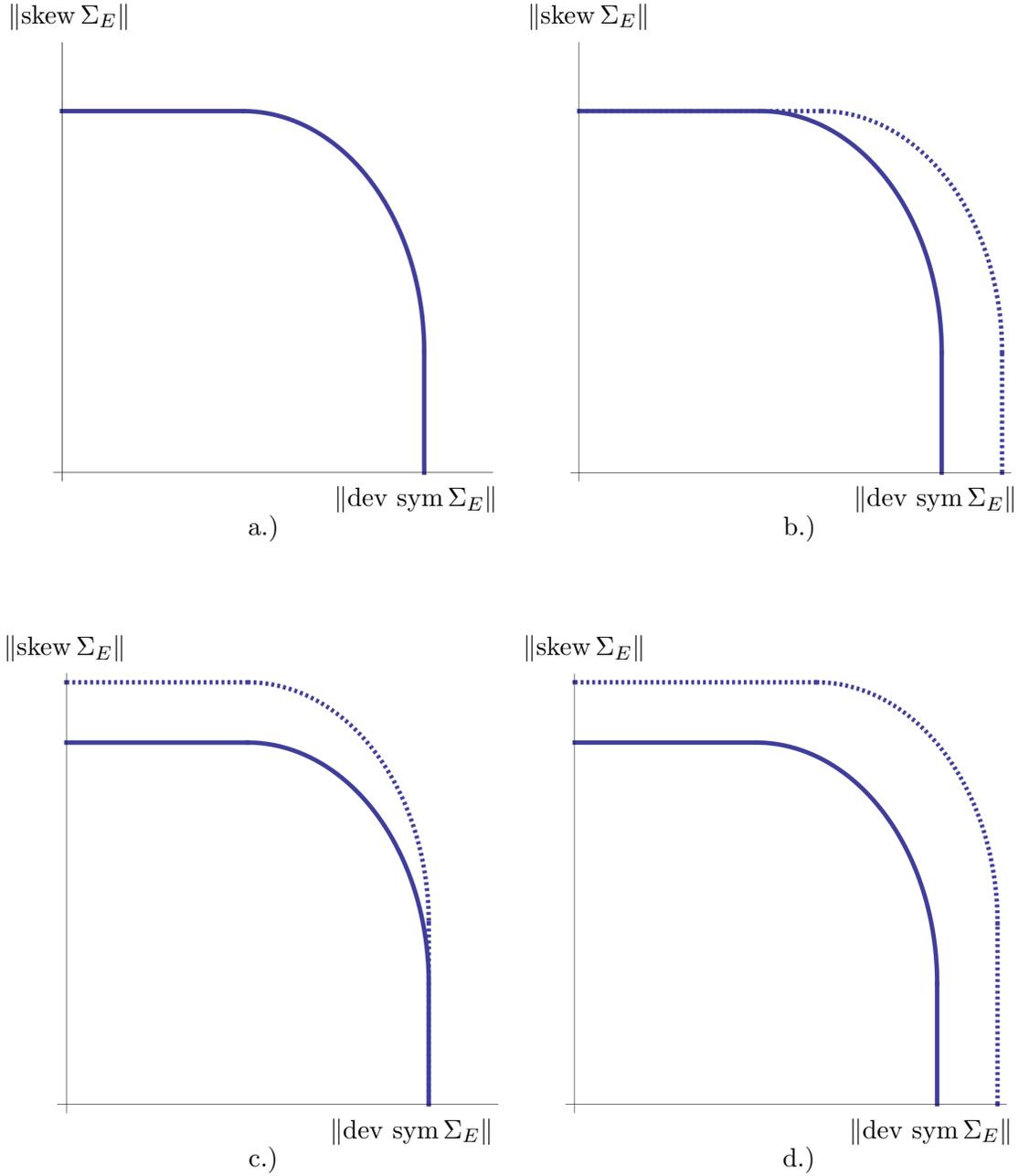
\begin{figure}[h!]
\centering
\vspace{0.5cm}
\includegraphics[width=0.4\textwidth]{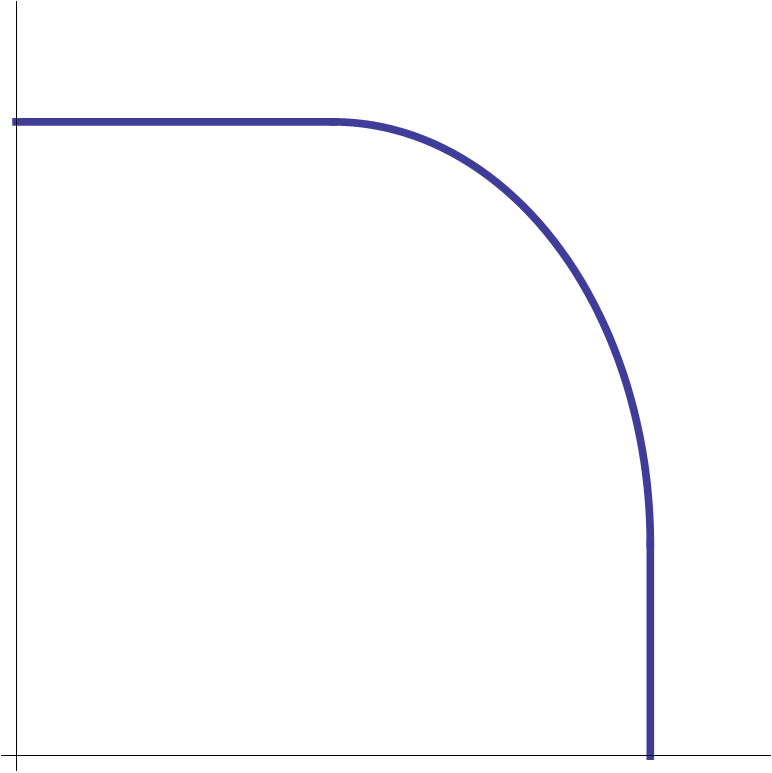} \hspace{1cm}\includegraphics[width=0.4\textwidth]{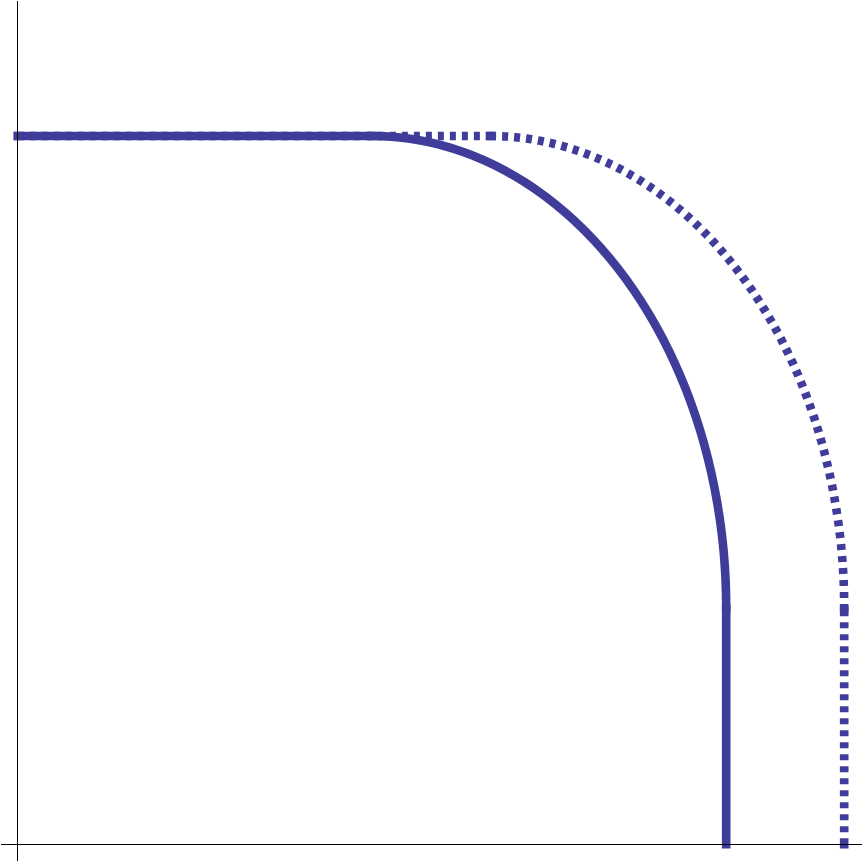}
\put(-280,-10){$\norm{\dev\,\sym\Sigma_E}$}
\put(-65,-10){$\norm{\dev\,\sym\Sigma_E}$}
\put(-415,190){$\norm{\skew\Sigma_E}$}
\put(-200,190){$\norm{\skew\Sigma_E}$}

a.) \hspace{7cm} b.)\\[2ex]
\vspace{1cm}

\includegraphics[width=0.4\textwidth]{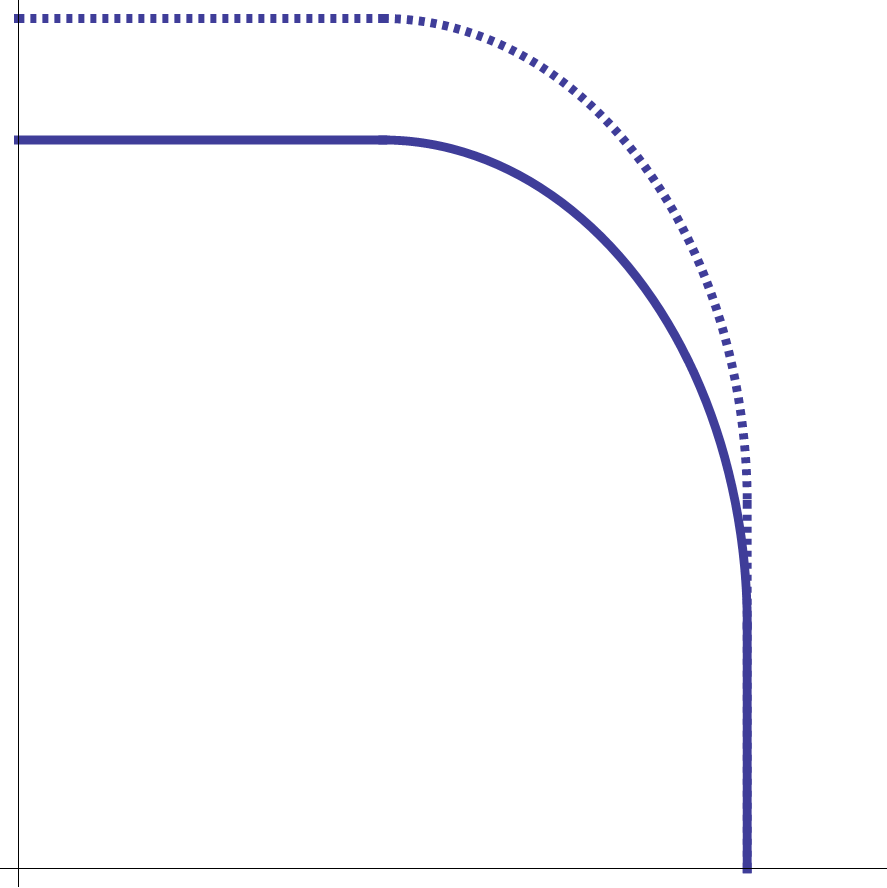}\hspace{1cm}\includegraphics[width=0.4\textwidth]{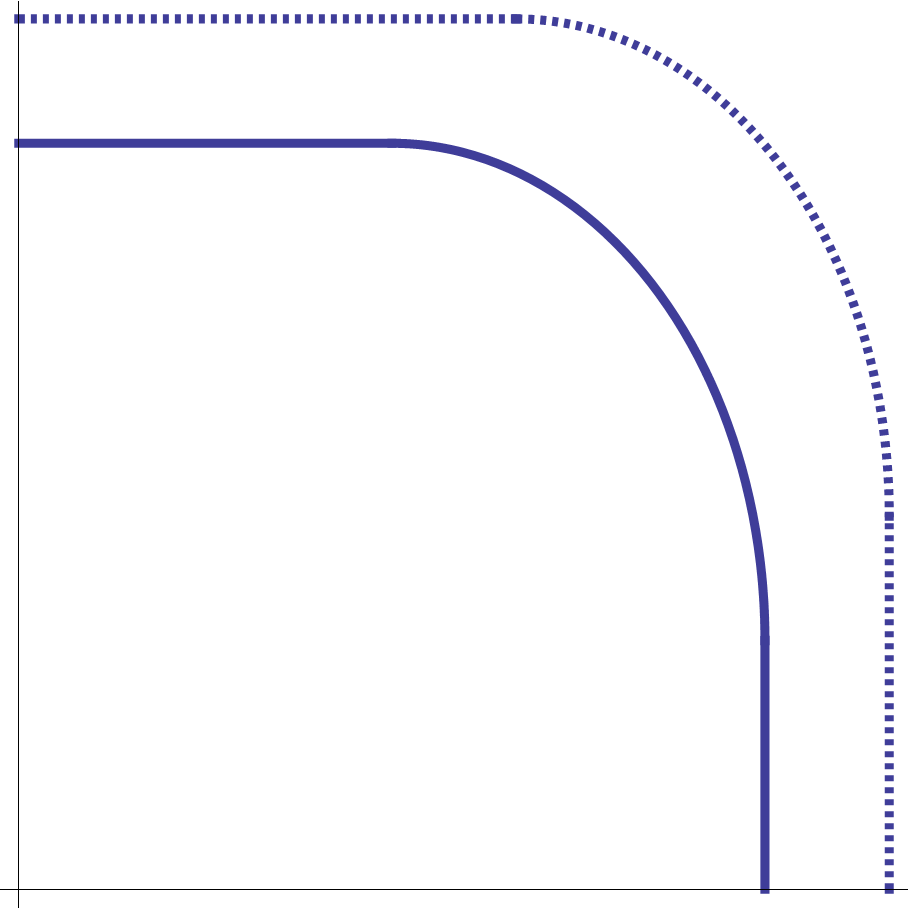}
\put(-280,-10){$\norm{\dev\,\sym\Sigma_E}$}
\put(-65,-10){$\norm{\dev\,\sym\Sigma_E}$}
\put(-415,190){$\norm{\skew\Sigma_E}$}
\put(-200,190){$\norm{\skew\Sigma_E}$}

c.) \hspace{7cm} d.)

\caption{Evolution of the yield surface by isotropic hardening: original yield surface depicted by solid line, evolved yield surface by dashed line. Different cases: a.) no evolution of yield surface, b.) expansion in direction of $\norm{\dev\,\sym\Sigma_E}$, c.) expansion in direction of $\norm{\skew\Sigma_E}$, d.) expansion in both directions.}
\label{fig:yield-surface}
\end{figure}

Let us briefly discuss the evolution of the yield surface. The hardening behavior will depend on the values of the moduli $\alpha_1, \alpha_2$, and on the location of the generalized stress state $(\Sigma_E,g_1,g_2)$ on the yield surface. There are four different possibilities, displayed in Fig. \ref{fig:yield-surface}. For clarity, we summarize the various cases in Table \ref{table:hardening}.

\begin{table}[h!]{\begin{center}

\begin{tabular}{|c|cccc|}
\hline
\rule[-1ex]{0pt}{2.5ex}   & $\alpha_1=0$, $\alpha_2=0$ & $\alpha_1>0$, $\alpha_2=0$ & $\alpha_1=0$, $\alpha_2>0$ & $\alpha_1>0$, $\alpha_2>0$ \\
\hline
\rule[-1ex]{0pt}{2.5ex} $(\Sigma_E,g_1,g_2) \in \mathcal{S}_1$ & a.) & b.) & a.) & b.) \\
\rule[-1ex]{0pt}{2.5ex} $(\Sigma_E,g_1,g_2) \in \mathcal{S}_2$ & a.) & a.) & c.) & c.) \\
\rule[-1ex]{0pt}{2.5ex} $(\Sigma_E,g_1,g_2) \in \mathcal{S}_3$ & a.) & b.) & c.) & d.) \\
\hline
\end{tabular}
\caption{Evolution of the yield surface, different hardening scenarios.}
\label{table:hardening}

\end{center}}
\end{table}

 Next, our goal is to present a strong and
a weak formulation of the model, followed by two existence results
for which there is an important distinction between the cases $\alpha_2>0$ and
$\alpha_2=0$ in the free-energy $\Psi$ in
(\ref{free-eng})-(\ref{free-eng-expr}).

\section{The complete mathematical formulation in the case $\alpha_2>0$}\label{alpha2-pos} In this section, we present the full description
of the model in the case $\alpha_2>0$ in the free-energy $\Psi$ in
(\ref{free-eng})-(\ref{free-eng-expr}) as well as a corresponding
existence result. The case $\alpha_2>0$ means that there is always
isotropic hardening in the spin-evolution equation. We recall that
the dissipation function $\Delta$ is given in (\ref{diss-funct1})
and the yield function in the case $\yieldzero>0$ and
$\yieldspino>0$ (see (\ref{S2})) is given by
\begin{equation}\label{yield-funct}\phi(\Sigma_p):=\left\{\begin{array}{ll}
\norm{\dev\sym\Sigma_E}+g_1-\yieldzero &\mbox{ \,on
}\mathcal{S}_1\\\\\norm{\skew\Sigma_E}+g_2-\yieldspino
&\mbox{ \,on }\mathcal{S}_2\\\\
\dsize\frac{(\norm{\dev\sym\Sigma_E}+g_1)^2}{\yieldzero^2}+
\dsize\frac{(\norm{\skew\Sigma_E}+g_2)^2}{\yieldspino^2}-1 &\mbox{
on }\mathcal{S}_3\,.\end{array}\right.\end{equation}

\subsection{The strong formulation}\label{strong} The strong formulation of the
model consists in finding:
\begin{itemize}\item[(i)] the displacement $u\in \SFH^1(0,T;
\SFH^1_0(\Omega,{\Gamma},\mathbb{R}^3))$,
\item[(ii)] the infinitesimal plastic distortion $\PL\in
\SFH^1(0,T;\SFL^2(\Omega, \sL(3)))$ with $$\Curl \PL\in \SFH^1(0,T;
\SFL^2(\Omega,\BBR^{3\times 3}))\quad\mbox{ and }\quad \Curl\Curl \PL
\in \SFH^1(0,T; \SFL^2(\Omega,\BBR^{3\times 3}))\,,$$\item[(iii)] the
internal isotropic hardening variables $\gamma_p,\,\,\omega_p\in\SFH^1(0,T;\SFL^2(\Omega))$
\end{itemize}
 such that the content of Table \ref{table:isohardspin} holds.

\begin{table}[h!]{\scriptsize\begin{center}
\begin{tabular}{|ll|}\hline &\qquad\\
 {\em Additive split of distortion:}& $\nabla u =e +\PL$,\quad $\bvarepsilon_e=\mbox{sym}\,e$,\quad $\bvarepsilon_p=\sym \PL$\\
{\em Equilibrium:} & $\mbox{Div}\,\sigma +f=0$ with
$\sigma=\C_{\mbox{\scriptsize iso}}\bvarepsilon_e$\\&\\ {\em Free
energy:} & $\frac12\langle\C_{\mbox{\scriptsize
iso}}\bvarepsilon_e,\bvarepsilon_e\rangle+\frac12\mu
\,L^2_c\,\norm{\mbox{Curl}\,\PL}^2+\frac12\mu\,\alpha_1\,|\gamma_p|^2+\frac12\mu\,\alpha_2\,|w_p|^2$\\&\\
{\em Yield condition:} & $\phi(\Sigma_p)=0$ with $\phi$ given
 in (\ref{yield-funct})\\&\\
 {\em where } & $\Sigma_p=(\Sigma_E,g_1,g_2)$,\,\quad $\Sigma_E:=\sigma+\Sigma^{\mbox{\scriptsize
lin}}_{\mbox{\scriptsize curl}}$,\,\quad $\Sigma^{\mbox{\scriptsize
lin}}_{\mbox{\scriptsize curl}}=-\mu \,L^2_c\,\Curl\Curl \PL$\\ &
$g_1=-\mu\,\alpha_1\gamma_p$,\,\quad \,$g_2=-\mu\,\alpha_2\,w_p$\\&
 \\{\em Dissipation inequality:} &
 $\dsize\int_\Omega[\langle\dev\sym\Sigma_E,\sym\dot{\PL}\rangle +\langle\skew\Sigma_E,\skew\dot{\PL}\ra +g_1\,\dot{\gamma_p}+g_2\,\dot{w}_p]
 \,dx\geq0$\\&\\
 {\em Dissipation function:} &$
\Delta(\dot{\Gamma}_p)$ is defined in (\ref{diss-funct1})\\&\\
 {\em Flow law in primal form:} &
 $\Sigma_p\in\partial \Delta(\dot{\Gamma}_p)$\\&\\
 {\em Flow law in dual form on $\mathcal{S}_1$:} &
$\left\{\begin{array}{ll}\sym\dot{p}\,\,=\lambda
\dsize\frac{\dev\sym\Sigma_E}{\norm{\dev\sym\Sigma_E}}
&\dot{\gamma}_p=\lambda=\norm{\sym\dot p}\quad
\\\skew\dot{p}=0 &\dot{\omega}_p=0\end{array}\right.$
 \\&\\
{\em Flow law in dual form on $\mathcal{S}_2$:} &
$\left\{\begin{array}{ll} \sym\dot{p}\,\,=0
&\,\,\quad\dot{\gamma}_p=0\\
\skew\dot{p}=\lambda\,\dsize\frac{\skew\Sigma_E}{\norm{\skew\Sigma_E}} &\quad\,\,\dot{\omega}_p=\lambda=\norm{\skew\dot p}\end{array}\right.$\\&\\

{\em Flow law in dual form on $\mathcal{S}_3$:}
&$\left\{\begin{array}{lll}\sym\dot{\PL}&=&\dsize\frac{2\lambda}{\yieldzero^2}\,(\norm{\dev\sym\Sigma_E}+g_1)
\dsize\frac{\dev\sym\Sigma_E}{\norm{\dev\sym\Sigma_E}}\quad
\\\\
\skew\dot{\PL}&=&\dsize\frac{2\lambda}{\yieldspino^2}\,(\norm{\skew\Sigma_E}+g_2)\dsize\frac{\skew\Sigma_E}{\norm{\skew\Sigma_E}}\\\\
\dot{\gamma}_p&=&\dsize\frac{2\lambda}{\yieldzero^2}\,(\norm{\dev\sym\Sigma_E}+g_1)=\norm{\sym\dot{\PL}}\\\\
\dot{\omega}_p&=&\dsize\frac{2\lambda}{\yieldspino^2}\,(\norm{\skew\Sigma_E}+g_2)=\norm{\skew\dot{\PL}}\\\\
\norm{\dot{\PL}}&=&\sqrt{\dot{\gamma}_p^2+\dot{\omega}_p^2}\,,\qquad
\lambda=\dsize\frac12\sqrt{\yieldzero^2\,\dot{\gamma}_p^2+\yieldspino^2\,\dot{\omega}_p^2}\,
\end{array}\right.$
\\&\\ {\em KKT conditions:}
&$\lambda\geq0$,\qquad $\phi(\Sigma_E,g_1,g_2)\leq0$,\qquad\quad
$\lambda\,\phi(\Sigma_E,g_1,g_2)=0$\\&\\
 {\em Boundary conditions for $\PL$:} & $\PL\times\,n=0$ on
 $\Gamma$,\,\, \quad $(\Curl \PL)\times\,n=0$ on $\partial\Omega\setminus\Gamma$\\
 {\em Function space for $\PL$:} & $\PL(t,\cdot)\in \mbox{H}(\mbox{Curl};\,\Omega,\,\BBR^{3\times 3})$\\
 \hline
\end{tabular}\caption{\footnotesize The model with isotropic hardening and plastic spin: the case $\alpha_2>0$. Because of the control of the $L^2$-norm of both
  isotropic hardening variables $\eta$ and $\beta$ and the constraints $\norm{\sym q}\leq\eta$ and $\norm{\skew q}\leq \beta$, the
  function space for the non-symmetric plastic distortion $p$ is easily seen to be $\mbox{H}(\mbox{Curl};\,\Omega,\,\BBR^{3\times 3})$.
  Note that the dual formulation of the flow rule needs a case distinction depending on the branches $\mathcal{S}_i$ of the
   yield surface while the primal formulation does not need it. Clearly, isotropic hardening for the plastic strain depends only on the accumulated equivalent plastic strain and isotropic hardening for the plastic
   rotation depends only on the accumulated equivalent plastic rotation. Therefore, there is no spin cross-hardening.
   In this flow rule in dual form we appreciate a Tresca like behaviour (see \cite{HEVALLER2005})
    in that we have to use a case distinction to determine on which part of the yield surface the evolution takes place.
}\label{table:isohardspin}
   \end{center}}\end{table}

\subsection{The weak formulation}\label{wf}

Assume that the problem in Section \ref{strong} has a solution
$w=(u,\Gamma_p)$ with $\Gamma_p=(\PL,\gamma_p,\omega_p)$.  Let $v\in H^1(\Omega,\mathbb{R}^3)$
with $v_{|\Gamma}=0$. Multiply the equilibrium equation with
$v-\dot{u}$ and integrate in space by parts and use the
symmetry of $\sigma$ and the elasticity relation to get
\begin{equation}\label{weak-eq1}
\int_{\Omega}\la\C_{\mbox{\scriptsize iso}}\sym(\nabla
u-\PL),\mbox{sym}(\nabla v-\nabla\dot{u})\ra\, dx=\int_\Omega
f(v-\dot{u})\,dx\, .
\end{equation}
 Now,
for any $q\in C^\infty(\overline{\Omega},\sL(3))$ such that
$q\times\,n=0$ on $\Gamma$ and any $\eta,\,\beta\in L^2(\Omega)$,
we integrate (\ref{flow-law-primal2}) over $\Omega$, integrate by
parts the term with Curl\,Curl using the boundary conditions
$$(q-\dot{\PL})\times\,n=0\mbox{ on }\Gamma,\qquad
\mbox{Curl}\,\PL\times\,n=0\mbox{ on }
\partial\Omega\setminus\Gamma$$ and get for $\overline{\Gamma}=(q,\eta,\beta)$
\begin{eqnarray}\label{flow-law-primal2-wk}
\nonumber &&\hskip-1truecm \int_\Omega\Delta(\overline{\Gamma})\,dx
-\int_\Omega\Delta(\dot{\Gamma}_p)\,dx -\int_\Omega
\la\C_{\mbox{\scriptsize iso}}(\sym\nabla u-\sym\PL),\sym q-\sym\dot{\PL}\ra\,dx\\
&&\quad+\mu\,L_c^2\int_\Omega\Bigl[\la\Curl
\PL,\Curl(q-\dot{\PL})\ra
+\mu\,\alpha_1\,\gamma_p\,(\eta-\dot{\gamma}_p) +\mu\,\alpha_2\,
\omega_p\,(\beta-\dot{\omega}_p)\,\Bigr]dx\geq0\,.\end{eqnarray}

Now adding up (\ref{weak-eq1}) and (\ref{flow-law-primal2-wk}) we
get the following weak formulation of the problem in Section
\ref{strong} in the form of a variational inequality:
\begin{eqnarray}\label{weak-form-iso}
&&\int_\Omega\la\C_{\mbox{\scriptsize iso}}(\sym\nabla
u-\sym\PL),(\sym\nabla v-\sym
q)-(\sym\nabla\dot{u}-\sym\dot{\PL})\ra\,dx\\
&&\nonumber\qquad+\,\,\mu\,L_c^2\int_\Omega\la\Curl
\PL,\Curl(q-\dot{\PL})\ra
+\,\int_\Omega\bigl[\mu\,\alpha_1\gamma_p(\eta-\dot{\gamma}_p)
+\mu\,\alpha_2\,\omega_p(\beta-\dot{\omega}_p)\bigr]\,dx\\
&&\nonumber\qquad\qquad\quad+\int_\Omega\Delta(\overline{\Gamma})\,dx
-\int_\Omega\Delta(\dot{\Gamma}_p)\,dx\,\geq\,\int_\Omega
f(v-\dot{u})\,dx\,.\nonumber\end{eqnarray}
\subsection{Existence result for the weak formulation}\label{existence}
To prove the existence result for the weak formulation
(\ref{weak-form-iso}), we closely follow the abstract machinery developed by
Han and Reddy in \cite{Han-ReddyBook} for mathematical problems in geometrically linear
classical plasticity and used for instance in \cite{DEMR1, REM, NCA,
EBONEFF, ENR2015} for models of gradient plasticity. Precisely, we will need the following theorem.
\begin{theorem}\label{han-reddy-theo}{\rm{\bf(\cite[Theorem 6.19]{Han-ReddyBook})}}\\
Let $\SFZ$ be a Hilbert space and let $\SFW$ be a nonempty closed convex cone of $\SFZ$. Consider the following problem: find $w\in\SFH^1([0,T]; \SFZ)$ with $w(0)=0$ such that for almost every $t\in[0,T]$, $\dot{w}(t)\in\SFW$ and
\begin{equation}\label{HR-ABS1}
\ba(w,z-\dot{w})+j(z)-j(\dot{w})\geq \langle
\ell,z-\dot{w}\rangle\mbox{ for every } z\in \SFW\,.\end{equation}
Assume that the following hold:\begin{itemize}\item[1.] the bilinear form $\ba$ is symmetric, continuous on $\SFZ$ and coercive on $\SFW$, i.e., there exist $C>0$ and $\alpha>0$ such that
\begin{equation}\label{HR-ABS2}
 \ba(w,z)\leq C\,\norm{w}_Z\,\norm{z}_Z\qquad\forall w,\,z\in\SFZ\quad\mbox{ and }\quad
\ba(z,z)\geq\alpha\,\norm{z}^2_Z\qquad\forall z\in \SFW\,;\end{equation}
\item[2.] $\ell\in \SFH^1([0,T];\SFZ')$  with $\ell(0)=0$.

\item[3.] the functional $j$ is non-negative, convex,  lower continuous  and positively $1$-homogeneous  $\SFW$, i.e., $j(sz)=|s|\,j(z)\quad\forall s\in\mathbb{R}\,,\quad\forall z\in\SFZ\,.$
\end{itemize} Then the problem (\ref{HR-ABS1}) has a solution $w\in \SFH^1([0,T];\SFZ)$.
\end{theorem}
\begin{remark}\label{initially-loaded}{\rm Theorem \ref{han-reddy-theo} above corresponds to models where the body is initially unlodaded that is, $\ell(0)=0$ and initially undeformed that is, $w(0)=0$. In the case of models with plastically pre-deformed bodies, Theorem  \ref{han-reddy-theo} is still valid with the initial condition $w(0)=w_0$ where $w_0\in W$ is the solution of the elliptic variational inequality
$$\ba(w_0,z-w_0)+j(z)-j(w_0)\geq \langle \ell(0),z-w_0\rangle\,\,\mbox{ for every }z\in\SFW\,.$$}\end{remark}\vskip.2truecm 
Now, in order to get an existence result for the weak formulation (\ref{weak-form-iso}) through the abstract result in  Theorem \ref{han-reddy-theo}, we write 
(\ref{weak-form-iso})  as (\ref{HR-ABS1}) with 
\begin{eqnarray}
\ba(w,z)&=&\int_{\Omega}\Bigl[\la\C_{\mbox{\scriptsize
iso}}(\sym\nabla u-\sym\PL),\sym\nabla
v-\sym q\ra+\mu\, L_c^2\la\Curl \PL,\Curl q\ra\\
\nonumber &&\hskip3truecm+\mu\,
\alpha_1\gamma_p\,\eta+\mu\,\alpha_2\,\omega_p\,\beta\Bigr]dx\,,
\label{bilin-iso-spin}
\\\nonumber\\
j(z)&=&\int_\Omega\Delta(q,\eta,\beta)\,dx\,,\label{functional-isospin}\\
\langle \ell,z\rangle&=&\int_\Omega
f\,v\,dx\,,\label{lin-form}\end{eqnarray} for
$w=(u,\PL,\gamma_p,\omega_p)$ and $z=(v,q,\eta,\beta)$ in
$\SFZ$\,.\\
The Hilbert space $\SFZ$ and the closed convex subset $\SFW$  will be 
constructed later  in such a way that the functionals $\ba$, $j$ and
$\ell$ satisfy the assumptions of Theorem \ref{han-reddy-theo}.

We let
\begin{eqnarray}
\nonumber\SFV&=&\mathsf{H}^1_0(\Omega,{\Gamma_{\mbox{\scriptsize
D}}},\mathbb{R}^3)=\{v\in \mathsf{H}^1(\Omega,\mathbb{R}^3)\,|\,
v_{|\Gamma}=0\}\,,\\
\nonumber \SFQ&=&\mbox{H}_0(\mbox{Curl};\,\Omega,\,\Gamma,\sL(3))\,,\label{space-p}\\
\nonumber \Lambda&=& L^2(\Omega)\,,\\
   \SFZ&=&\SFV\times \SFQ\times\Lambda^2\,,\label{product-space}\\
   \SFW&=&\bigl\{z=(v,q,\eta,\beta)\in\SFZ\,\,|\,\,\norm{\sym\,q}\leq\eta\mbox{
   and }\norm{\skew\,q}\leq\beta\bigr\}\,,\label{set-W}
\end{eqnarray} and define the norms
\begin{eqnarray}
\nonumber&&  \norm{v}_V:=\norm{\nabla v}_{L^2},\label{norm-V}\quad\qquad
\norm{q}_{Q}:=\norm{q}_{\mbox{\scriptsize H}(\mbox{\scriptsize
Curl};\Omega)},\label{norm-Q}\\
&&\norm{z}^2_{Z}:=\norm{v}^2_{V} +\norm{q}^2_{Q}
+\norm{\eta}^2_{L^2}+\norm{\beta}^2_{L^2}\quad\mbox{ for }
z=(v,q,\eta,\beta)\in \SFZ\,. \label{norm-Z}
\end{eqnarray}
We will assume that the body is initially unloaded and undeformed and this corresponds to assuming that $f(x,0)=0$ for almost all $x\in\Omega$ with homogeneous initial conditions. We then get the following existence result for the weak formulation (\ref{weak-form-iso}).
\begin{theorem}\label{existence-wf-iso} Under the choices of the Hilbert space $\SFZ$ and the closed convex cone $\SFW$ in (\ref{space-p})-(\ref{set-W}) with the norms in (\ref{norm-Z}) and the functionals $\ba$, $j$ and
$\ell$ in (\ref{bilin-iso-spin})-(\ref{lin-form}), the weak formulation (\ref{weak-form-iso}) when written as the variational inequality of the second kind (\ref{HR-ABS1}) has a solution $w=
(u,\PL,\gamma_p,\omega_p) \in \SFH^1([0,T];\SFZ)$ with $\dot{w}\in\SFL^2([0,T];\SFW)$.\end{theorem}
 {\bf Proof:} The functionals $j$ and $\ell$ trivially satisfy the asumptions of Theorem \ref{han-reddy-theo}. So,  the key issue here is the
coercivity of the bilinear form $\ba$ on the set $\SFW$. 
 Let therefore $z=(v,q,\eta,\beta)\in \SFW$.
\begin{eqnarray}
\nonumber\ba(z,z)& \geq&  m_0\norm{\sym\,\nabla v-\sym q}^2_2\mbox{
(from (\ref{ellipticityC}))}\\ \nonumber &&\hskip3truecm  +\,\mu\,
L_c^2\norm{\Curl q}_2^2 +\mu\, \alpha_1\norm{\eta}_2^2+\mu\,\alpha_2\norm{\beta}_2^2 \\
\nonumber&=&m_0\left[\norm{\sym \,\nabla v}^2_2 +\norm{\sym q}^2_2
-2\la\sym
\,\nabla v, \sym \PL\ra\right]\\
\nonumber&&\qquad\qquad+\,\mu L_c^2\,\norm{\Curl q}_2^2
+\mu\, \alpha_1\norm{\eta}_2^2+\mu\,\alpha_2\norm{\beta}_2^2\\
\nonumber&\geq &m_0\left[\norm{\sym\,\nabla v}^2_2 +\norm{\sym
q}^2_2 -\theta\norm{\sym \,\nabla
v}_2^2-\frac1\theta\norm{\sym q}_2^2\right]\mbox{ (Young's inequality)}\\
\nonumber&& \qquad\qquad+\,\mu \,L_c^2\norm{\Curl q}_2^2+\frac12\mu\, \alpha_1\norm{\eta}_2^2+\frac12\mu\,\alpha_2\norm{\beta}_2^2 \\
\nonumber&&\qquad\qquad +\,\frac12\mu\, \alpha_1\norm{\sym q}_2^2
+\frac12\mu\,
\alpha_2\norm{\skew q}_2^2\mbox{ \,\,(using $\norm{\sym q}\leq\eta$,\,\, $\norm{\skew q}\leq\beta$)} \\
\nonumber &=& m_0(1-\theta)\norm{\sym \,\nabla
v}^2_2+\left[m_0\Bigl(1-\frac1\theta\Bigr)+\frac12\,\mu\,
\alpha_1\right]\norm{\sym q}_2^2+\frac12\mu\, \alpha_2\norm{\skew
q}_2^2\\
&&\qquad\qquad +\,\mu \,L_c^2\norm{\Curl q}_2^2 +
\frac12\,\mu\,\alpha_1\norm{\eta}_2^2+\frac12\,\mu\,\alpha_2\norm{\beta}_2^2.\label{coercive-a}\end{eqnarray}

So, choosing $\theta$ such that $\displaystyle\frac{m_0}{m_0+\frac
12\,\mu\,\alpha_1}< \theta<1,$ and using Korn's first inequality (see e.g. \cite{NEFFKORN2002}),
there exists some positive constant
$C(m_0,\mu,\alpha_1,\alpha_2,L_c,\Omega)>0$ such that
$$a(z,z)\geq C\left[\norm{v}_V^2+\norm{q}^2_{\mbox{\scriptsize H}(\mbox{\scriptsize
Curl};\Omega)}+\norm{\eta}^2_2+\norm{\beta}_2^2
\right]=C\norm{z}^2_{Z}\quad\forall z=(v,q,\eta,\beta)\in \SFW\,.$$ Hence, we get from Theorem \ref{han-reddy-theo} the existence of a solution  for the weak formulation (\ref{weak-form-iso})  with $\alpha_2>0$ (and $\alpha_1>0$).\hfill\qed
\subsection{Uniqueness of the strong solution}\label{uniqueness}

If in the geometrically linear classical plasticity model with isotropic hardening, the uniqueness of the weak solution is
 obtained from the formulation in a variational inequality (see \cite[Theorem 7.3]{Han-ReddyBook})
 the uniqueness of the weak solution in the context of gradient plasticity with
 isotropic hardening has not yet been completely established.  
 However, in some particular cases,
  the uniqueness has been obtained provided weak solutions are regular enough (see e.g. \cite[pp.210-212]{Han-ReddyBook}). 
  
  The diffculty here is that the coerciviy of the bilinear form $\ba$, which is key to get the uniqueness of the solution of the weak formulation, is only obtained on the closed convex cone $\SFW$ and not on the  entire space $\SFZ$. Therefore, one cannot use the standard argument of involving the difference $w_1-w_2$ of two solutions $w_1$ and $w_2$ and getting the uniqueness, since that difference does not always belong to the closed convex cone $\SFW$. We recall that, for Prager-type linear kinematical hardening, the uniqueness
   of strong solutions in infinitesimal  perfect
   gradient plasticity was established in \cite{NEFF-IUTAM08}. In our context,  we will prove in the next theorem      that requiring
 $\Curl\Curl p\in \SFH^1([(0,T];\SFL^2(\Omega,\mathbb{R}^{3\times 3}))$ is enough to guarantee the uniqueness of the strong solution.
 \begin{theorem}\label{uniqueness-strong}
Let  $w=(u,\Gamma_p)$  be a solution of the weak formulation (\ref{weak-form-iso}) obtained in Theorem \ref{existence-wf-iso} with $\Gamma_p=(p,\gamma_p,\omega_p)$. If $p$ satisfies $\Curl\Curl p\in \SFH^1([0,T];\SFL^2(\Omega,\mathbb{R}^{3\times 3}))$, then $w$ is the unqiue strong solution, i.e., the unique solution of the strong formulation in Section \ref{strong}.
 \end{theorem}
{\bf Proof:}  In fact, we first notice that if $w=(u,\Gamma_p)\in \SFW$ is a solution of
(\ref{weak-form-iso}) with $\Gamma_p=(p,\gamma_p,\omega_p)$ and $\Curl\Curl\,{p}\in
L^2(\Omega,\mathbb{R}^{3\times 3})$, then choosing appropriately
test functions and integrating by parts, we easily get that $w=(u,\Gamma_p)$
 satisfies the equilibrium
equation (\ref{equil}) on the one hand and $\Gamma_p$ satisfies the flow rule
 in dual form
 \begin{equation}\label{normality}
\la\dot{\Gamma}_p,\overline{\Sigma}-\Sigma_p\ra\leq0\qquad\forall\,\overline{\Sigma}\,\end{equation}
on the other hand.

Let us now consider two solutions $w_i=(u_i,\Gamma_{p_i})$ \,
$i=1,2$ of (\ref{weak-form-iso}) with 
$\Gamma_{p_i}=(p_i,\gamma_{p_{\,i}},\omega_{p_{\,i}})$,  satisfying the same
 initial conditions and let $\Sigma_{p_{\,i}}=(\Sigma_{E_{\,i}},g_{1_{\,i}},g_{2_{\,i}})$ be the corresponding stresses. That is,
 \begin{equation}
 \Sigma_{E_{\,i}}=\sigma_i-\mu\,L^2_c\Curl\Curl p_i\,,\qquad g_{1_{\,i}}=-\mu\,\alpha_1\gamma_{p_{\,i}}\,,\qquad
 g_{2_{\,i}}=-\mu\,\alpha_2\,\omega_{p_{\,i}}\,.\end{equation}
Hence, $\Gamma_{p_{\,i}}$ and $\Sigma_{p_{\,i}}$ satisfy
$(\ref{normality})$, that is,
\begin{equation}\label{normality1-2}
\la{\dot\Gamma}_{p_1},\overline{\Sigma}-\Sigma_{p_{\,1}}\ra\leq0\quad\mbox{and}\quad
\la\dot{\Gamma}_{p_2},\overline{\Sigma}-\Sigma_{p_{\,2}}\ra\leq0\qquad\forall\,\overline{\Sigma}\,\end{equation}
Now choose $\overline{\Sigma}=\Sigma_{p_{\,2}}$ in
(\ref{normality1-2})$_1$ and $\overline{\Sigma}=\Sigma_{p_{\,1}}$ in
(\ref{normality1-2})$_2$ and add up to get
\begin{equation}\label{diff-norm1-2}
\la\Sigma_{p_{\,2}}-\Sigma_{p_{\,1}},\dot{\Gamma}_{p_{\,1}}-\dot{\Gamma}_{p_{\,2}}\ra\leq0\,.\end{equation}
That is
\begin{eqnarray}\label{diff-normality1}
&&\nonumber \hskip-1truecm\la\sigma_2-\sigma_1,  \dot{p}_1-\dot{p}_2\ra +\mu\,L_c^2\la\Curl\Curl(p_2-p_1),\dot{p}_2-\dot{p}_1\ra
\\
&&\hskip1truecm
+\,\mu\,\alpha_1\,(\gamma_{p_{\,2}}-\gamma_{p_{\,1}})(\dot{\gamma}_{p_{\,2}}-\dot{\gamma}_{p_{\,1}})+
\mu\,\alpha_2\,(\omega_{p_{\,2}}-\omega_{p_{\,1}})(\dot{\omega}_{p_{\,2}}-\dot{\omega}_{p_{\,1}})\leq0\,.\end{eqnarray}
Since $\sigma$ is symmetric, the latter is equivalent to
\begin{eqnarray}\label{diff-normality11}
&&\nonumber \hskip-1truecm\la\sigma_2-\sigma_1,
\sym(\dot{p}_1-\dot{p}_2)\ra
+\mu\,L_c^2\la\Curl\Curl(p_2-p_1),\dot{p}_2-\dot{p}_1\ra
\\
&&\hskip.5truecm
+\,\mu\,\alpha_1\,(\gamma_{p_{\,2}}-\gamma_{p_{\,1}})(\dot{\gamma}_{p_{\,2}}-\dot{\gamma}_{p_{\,1}})+
\mu\,\alpha_2\,(\omega_{p_{\,2}}-\omega_{p_{\,1}})(\dot{\omega}_{p_{\,2}}-\dot{\omega}_{p_{\,1}})\leq0\,.\end{eqnarray}
 Now, substitute $\sym p_i=\sym(\nabla u_i)-\C^{-1}\sigma_i$
obtained from the elasticity relation, into equation (\ref{diff-normality11})
and get
\begin{eqnarray}\label{diff-normality2}
&&\nonumber \hskip-1truecm\la\sigma_2-\sigma_1,
\C^{-1}(\dot{\sigma}_2-\dot{\sigma}_1)\ra
+\mu\,L_c^2\la\Curl\Curl(p_2-p_1),\dot{p}_2-\dot{p}_1\ra
 +\,\mu\,\alpha_1\,(\gamma_{p_{\,2}}-\gamma_{p_{\,1}})(\dot{\gamma}_{p_{\,2}}-\dot{\gamma}_{p_{\,1}})\\
&&\hskip1truecm
+\,\mu\,\alpha_2\,(\omega_{p_{\,2}}-\omega_{p_{\,1}})(\dot{\omega}_{p_{\,2}}-\dot{\omega}_{p_{\,1}})
\leq\la\sigma_1-\sigma_2,\sym(\nabla\dot{ u}_1)-\sym(\nabla
\dot{u}_2)\ra\,.\end{eqnarray} Now, notice that from the equilibrium
equation we get
$$\int_\Omega\la\sigma_1-\sigma_2,\sym(\nabla\dot{u}_1)-\sym(\nabla\dot{u}_2)\ra\,dx=0\,.$$
Hence, for a.e. $t\in(0,T)$, integrate (\ref{diff-normality2}) over $\Omega\times(0,t)$ then after integrating the term with $\Curl\Curl$ by parts, we get
\begin{eqnarray}\label{diff-normality3}
&&\nonumber\hskip-1truecm\int_0^t\frac{d}{ds}\Bigl[\left.\norm{\C^{-1/2}(\sigma_2(s)-\sigma_1(s))}^2_2+\mu\,L^2_c\norm{\Curl(p_1(s)-p_2(s))}^2_2\right.\\
&&\hskip2truecm
\left.+\,\mu\,\alpha_1\,\norm{\gamma_{p_{\,2}}(s)-\gamma_{p_{\,1}}(s)}_2^2
+\mu\,\alpha_2\,\norm{\omega_{p_{\,1}}(s)-\omega_{p_{\,2}}(s)}_2^2\right.\Bigr]\,ds\leq0\,.
\end{eqnarray} Therefore, we obtain
\begin{eqnarray}\label{unique-ineq}
\nonumber
&&\hskip-2truecm\norm{(\C^{-1})^{1/2}(\sigma_2-\sigma_1)}^2_2+\mu\,L^2_c\,\norm{\Curl(p_1-p_2)}^2_2\\
&&\qquad\qquad\quad+\,\,\mu\,\alpha_1\,\norm{\gamma_{p_{\,2}}-\gamma_{p_{\,1}}}_2^2
+\mu\,\alpha_2\,\norm{\omega_{p_{\,2}}-\omega_{p_{\,1}}}_2^2=0\,,
\end{eqnarray} from which we get $\sigma_1=\sigma_2$, \,$\Curl p_1=\Curl p_2$,\, $\gamma_{p_{\,1}}=\gamma_{p_{\,2}}$,\,
$\omega_{p_{\,2}}=\omega_{p_{\,1}}$ and hence,
$\Sigma_{E_{\,1}}=\Sigma_{E_{\,2}}$.

Now, let us prove that $p_1=p_2$. In fact, from the definition of
the normal cone it follows that $\dot{p}_i=0$ and
$\dot{\gamma}_{p_{\,i}}=\dot{\omega}_{p_{\,i}}=0$ inside the elastic
domain $\mathcal{E}$, which from the initial conditions imply that
$p_i=0$ inside $\mathcal{E}$. Now, looking at the flow rule in dual
form in Table \ref{table:isohardspin}, we easily obtain that
$\sym\dot{p_1}=\sym\dot{p}_2$ and $\skew\dot{p}_1=\skew\dot{p}_2$ on
each surface $S_k$. Therefore, $\dot{p}_1=\dot{p}_2$ which implies
that $p_1=p_2$ from the initial conditions.

 In order to show that $u_1=u_2$, we use $\sym(\nabla u_i)=\C^{-1}\sigma_i+\sym p_i$ obtained from the elasticity relation and get
 $$\sym(\nabla (u_1-u_2))=\C^{-1}(\sigma_1-\sigma_2)+\sym(p_1-p_2)=0\,,$$
 and hence, from the first Korn's inequality (see e.g. \cite{NEFFKORN2002}), we get $\nabla(u_1-u_2)=0$ which implies that $u_1=u_2$. Therefore, we finally obtain
 $$u_1=u_2\,,\qquad \sigma_1=\sigma_2\,,\qquad p_1=p_2\,,\qquad \gamma_{p_{\,1}}=\gamma_{p_{\,2}}\,,
 \qquad\omega_{p_{\,1}}=\omega_{p_{\,2}}\,,$$ and thus the uniqueness
  of a strong solution to the mathematical problem describing our model of rate-independent geometrically linear gradient plasticity with
  isotropic hardening and plastic spin in the case $\alpha_2>0$, where there is always isotropic hardening in the spin-evolution equation.\qquad$\mbox{ }$\hfill\qed
\subsection{Perfect gradient plasticity with
spin}\label{perfect-graplast} Inspection of the uniqueness proof for
strong solutions in Section \ref{uniqueness}  shows that in the case
with zero isotropic plastic strain and spin hardening, and the homogeneous boundary conditions $u|_\Gamma=0$ and $p\times n|_{\Gamma}=0$,  elastic
stresses $\sigma$, elastic strains $\bvarepsilon_e=\sym e$ and furthermore elastic
distortions $e=\nabla u-p$ are unique. The uniqueness with respect
to elastic distortions uses again the new Korn's inequality for
incompatible tensor fields \cite{NPW2014} since $e\times n|_\Gamma=0$. In this case, the extra inclusion of the
spin and the dislocation density tensor allow to improve uniqueness
from elastic strains to elastic distortions. Notice that, in the context of crystal gradient plasticity, non-uniqueness has been shown in \cite{Bardegiacomini2008} for the case of nonhomogeneous displacement boundary conditions, focussing on the simple shear of a constrained strip endowed with
multiple slip systems. The same type of non-unqiueness results have been obtained also in \cite{GURT2004,POHPEERLINGS2016,PANTEBARDE2018}.
\section{The complete mathematical formulation in the no-spin-hardening
case}\label{alpha2-zero} Here we set $\alpha_2=0$ in the free-energy
$\Psi$ in (\ref{free-eng})-(\ref{free-eng-expr}). The case $\alpha_2=0$ means that there is no isotropic hardening in the spin-evolution.
 At present we believe that it is this case which deserves special attention, since in this model we extend classical plasticity in the
  weakest possible way to depend on plastic spin. Notably, we do not incur additional spin-hardening.
 The dissipation
function $\Delta$ is still the same given in (\ref{diss-funct1}) and
the yield function is given in (\ref{yield-funct}). Also, in
this model the influence of the SSD's and GND's on plastic flow is
neatly separated: the SSD-distribution influences only isotropic
hardening through the classical mechanism and the GND-distribution
determines the nonlocal kinematic hardening.
\subsection{The strong formulation of the model}\label{sf-spinfree}
The strong formulation in the no-spin-hardening case is obtained  exactly as in Section \ref{strong}.  For the clarity of exposition, we chose to present here the whole formulation summarized in Table  \ref{table:isohardspin-free} below, instead of just pointing out the differences w.r.t. Table  \ref{table:isohardspin}.
\begin{table}[h!]\footnotesize\begin{center}
\begin{tabular}{|ll|}\hline &\qquad\\
 {\em Additive split of distortion:}& $\nabla u =e +\PL$,\quad $\bvarepsilon_e=\mbox{sym}\,e$,\quad $\bvarepsilon_p=\sym \PL$\\
{\em Equilibrium:} & $\mbox{Div}\,\sigma +f=0$ with
$\sigma=\C_{\mbox{\scriptsize iso}}\bvarepsilon_e$\\&\\ {\em Free
energy:} & $\frac12\langle\C_{\mbox{\scriptsize
iso}}\bvarepsilon_e,\bvarepsilon_e\rangle+\frac12\mu
\,L^2_c\,\norm{\mbox{Curl}\,\PL}^2+\frac12\mu\,\alpha_1\,|\gamma_p|^2$\\&\\
{\em Yield condition:} & $\phi(\Sigma_p)=0$ with $\phi$ given
 in (\ref{yield-funct})\\
 {\em where } & $\Sigma_p=(\Sigma_E,g_1,g_2)$,\,\qquad $\Sigma_E:=\sigma+\Sigma^{\mbox{\scriptsize
lin}}_{\mbox{\scriptsize curl}}$,\,\qquad $\Sigma^{\mbox{\scriptsize
lin}}_{\mbox{\scriptsize curl}}=-\mu \,L^2_c\,\Curl\Curl \PL$\\ &
$g_1=-\mu\,\alpha_1\gamma_p$,\,\qquad $g_2=0$\\&
 \\{\em Dissipation inequality:} &
 $\dsize\int_\Omega [\langle\dev\sym\Sigma_E,\sym\dot{\PL}\rangle +\la\skew\Sigma_E,\skew{\dot p}\ra +g_1\,\dot{\gamma}_p]
 \,dx\geq0$\\
 {\em Dissipation function:} &$
\Delta(\dot{\Gamma}_p)$ is defined in (\ref{diss-funct1})\\
 {\em Flow law in primal form:} &
 $\Sigma_p\in\partial \Delta(\dot{\Gamma}_p)$\\&\\
 {\em Flow law in dual form on $\mathcal{S}_1$:} &
$\left\{\begin{array}{ll}\sym\dot{p}\,\,=\lambda\,\dsize\frac{\dev\sym\Sigma_E}{\norm{\dev\sym\Sigma_E}},&\dot{\gamma}_p=\lambda=\norm{\sym\dot
p}\\ \skew\dot{p}=0,&\dot{\omega}_p=0\end{array}\right. $
 \\&\\

{\em Flow law in dual form on $\mathcal{S}_2$:} &
$\left\{\begin{array}{ll}\sym\dot{p}\,\,=0, & \quad\,\,\dot{\gamma}_p=0\\
\skew{\dot
p}=\lambda\,\dsize\frac{\skew\Sigma_E}{\norm{\skew\Sigma_E}}, &\quad\,\,\dot{\omega}_p=\lambda=\norm{\skew\dot{p}}\end{array}\right.$\\&\\

{\em Flow law in dual form on $\mathcal{S}_3$:}
&$\left\{\begin{array}{lll}\sym\dot{\PL}&=&\dsize\frac{2\lambda}{\yieldzero^2}\,(\norm{\dev\,\sym\,\Sigma_E}+g_1)
\dsize\frac{\dev\sym\Sigma_E}{\norm{\dev\sym\Sigma_E}}\quad
\\\\
\skew\dot{\PL}&=&\dsize\frac{2\lambda}{\yieldspino^2}\skew\Sigma_E\\\\
\dot{\gamma}_p&=&\dsize\frac{2\lambda}{\yieldzero^2}\,(\norm{\dev\,\sym\,\Sigma_E}+g_1)=\norm{\sym\dot{\PL}}\\\\
\dot{\omega}_p&=&\dsize\frac{2\lambda}{\yieldspino^2}\,\norm{\skew\,\Sigma_E}=\norm{\skew\dot{\PL}}\\\\
\norm{\dot{\PL}}&=&\dsize\sqrt{\dot{\gamma}_p^2+\dot{\omega}_p^2}\,,\qquad
\lambda=\dsize\frac12\sqrt{\yieldzero^2\,\dot{\gamma}_p^2+\yieldspino^2\,\dot{\omega}_p^2}\,
\end{array}\right.$
\\&\\ {\em KKT conditions:}
&$\lambda\geq0$,\qquad $\phi(\Sigma_E,g_1,0)\leq0$,\qquad
$\lambda\,\phi(\Sigma_E,g_1,0)=0$\\&\\
 {\em Boundary conditions for $\PL$:} & $\PL\times\,n=0$ on
 $\Gamma$,\,\,\quad $(\Curl \PL)\times\,n=0$ on
 $\partial\Omega\setminus\Gamma$\\
 {\em Function space for $\PL$:} & $\PL(t,\cdot)\in \mbox{H}(\mbox{Curl};\,\Omega,\,\BBR^{3\times 3})$\\
 \hline
\end{tabular}\caption{\footnotesize The model with isotropic hardening only in the plastic strain-evolution (the case $\alpha_2=0$). Notice that the boundary
 conditions on $p$ necessitates at least $p\in \mbox{H}(\Curl;\,\Omega,\,\BBR^{3\times3})$, which is not guaranteed looking at the
 free-energy and the dissipation function. However, this will be obtained from a new Korn's type inequality  for incompatible tensor fields derived by Neff et al. in
  \cite{NPW2011-1, NPW2012-1, NPW2012-2,
NPW2014}\,.}\label{table:isohardspin-free}\end{center}\end{table}
\subsection{The weak formulation of the model}\label{wf-spinfree}
Also, following Section \ref{wf}, we derive the weak formulation of the model in the no-spin-hardening case as the variational inequality 
\begin{eqnarray}\label{weak-form-iso2}
&&\nonumber \hskip-0.7truecm \int_\Omega\la\C_{\mbox{\scriptsize iso}}(\sym\nabla
u-\sym\PL),(\sym\nabla v-\sym
q)-(\sym\nabla\dot{u}-\sym\dot{\PL})\ra\,dx +\alpha_1\,\mu\int_\Omega\gamma_p(\eta-\dot{\gamma}_p)\,dx\\
&& +\,\,\mu\,L_c^2\int_\Omega\la\Curl
\PL,\Curl(q-\dot{\PL})\ra\,dx+\int_\Omega\Delta(\overline{\Gamma})\,dx
-\int_\Omega\Delta(\dot{\Gamma}_p)\,dx\,\geq\,\int_\Omega
f(v-\dot{u})\,dx\,\,.\end{eqnarray}
\subsection{Existence result in the no-spin-hardening case}\label{existence-nospin-hardening}
As in Section \ref{existence}, the existence result in the no-spin-hardening is also obtained through the abstract result in Theorem \ref{han-reddy-theo}.  The functionals $j$ and $\ell$ remain as in (\ref{functional-isospin}) and (\ref{lin-form}) respectively, the bilinear form $\ba$ in this case is defined as
\begin{eqnarray}
\nonumber \ba(w,z)&=&\int_{\Omega}\Bigl[\la\C_{\mbox{\scriptsize
iso}}(\sym\nabla u-\sym\PL),\sym\nabla v-\sym q)\ra\\
&&\hskip4truecm+\mu\, L_c^2\la\Curl \PL,\Curl q\ra+\mu\,
\alpha_1\gamma_p\,\eta\Bigr]dx\,, \label{bilin-iso-special}
\\\nonumber\\
&&\nonumber\qquad\forall w=(u,\PL,\gamma_p,\omega_p),\,\,z=(v,q,\eta,\beta) \in\SFZ\,.\end{eqnarray}  The existence result for the weak formulation  in the no-spin-hardening case is obtained in the following theorem.
\begin{theorem}\label{existence-nospin-hardeding} The weak formulation (\ref{weak-form-iso2}) in the no-spin-hardening case (i.e., $\alpha_2=0$) written as: find 
$w=(u,\PL,\gamma_p,\omega_p)\in \SFH^1(0,T;\SFZ)$ such that $w(0)=0$
and $\dot{w}(t)\in \SFW$ for a.e. $t\in [0,T]$
\begin{equation}\label{wf-special}
\ba(\dot{w},z-w)+j(z)-j(\dot{w})\geq \langle
\ell,z-\dot{w}\rangle\mbox{ for every } z\in \SFZ\mbox{ and for a.e.
}t\in[0,T]\,,\end{equation} with $\ba$ defined in (\ref{bilin-iso-special}) and $j$ and $\ell$ defined in  (\ref{functional-isospin})-(\ref{lin-form}), has a solution for some suitable Hilbert space $\SFZ$ and some closed 
convex cone $\SFW$ in $\SFZ$.\end{theorem}
{\bf Proof:}
First of all, notice that since the bilinear form $\ba$ does not
contain explicitly the variable $\beta$, it is impossible to derive
the coercivity of the bilinear form in any normed space in all the
variables $v$, $q$, $\eta$ and $\beta$. Therefore, we are not in a position to apply directly the abstract result in Theorem \ref{han-reddy-theo}.  The new solution
strategy here for the existence result is to first find $u$, $p$ and
$\gamma_p$, and construct $\omega_p$ a posteriori. To this end, we
define

\begin{equation}\label{diss-funct-new}
\Delta_0(q,\eta):=\left\{\begin{array}{ll}\hskip-.2truecm
D(\norm{\sym q},\norm{\skew q})
& \mbox{if }\norm{\sym q}\leq\eta\,,\\\\
\hskip-.2truecm\infty &
\mbox{otherwise}\,,\end{array}\right.\end{equation} where we recall
that
\begin{equation}\label{pre-diss0}
D(s,t):=\sqrt{\yieldzero^2\,s^2+\yieldspino^2\,t^2}\,.
\end{equation}
We then reformulate the problem as follows: find
$w=(u,\PL,\gamma_p)\in \SFH^1(0,T;\SFZ)$ such that $w(0)=0$,\,
$\dot{w}(t)\in \SFW$ for a.e. $t\in [0,T]$
\begin{equation}\label{wf-new1}
\ba(\dot{w},z-w)+j_0(z)-j_0(\dot{w})\geq \langle
\ell,z-\dot{w}\rangle\mbox{ for every } z=(v,q,\eta)\in \SFW\mbox{
and for a.e. }t\in[0,T]\,,\end{equation} where we let
\begin{equation}\label{new-j}
j_0(z):=\left\{\begin{array}{ll}
\dsize\int_\Omega\Delta_0(q,\eta)\,dx
&\mbox{ if }z=(v,q,\eta)\in \SFW\\\\
\infty &\mbox{ otherwise}\,,\end{array}\right.\end{equation}
\begin{eqnarray}
\SFZ&=&\SFV\times \SFQ\times\Lambda\,,\label{product-space-new}\\
   \SFW&=&\bigl\{z=(v,q,\eta)\in\SFZ\,\,|\,\,\norm{\sym\,q}\leq\eta\,\,\mbox{ a.e. in }\,\Omega\bigr\}\,,\label{set-W-new}\\
\nonumber\SFV&=&\mathsf{H}^1_0(\Omega,{\Gamma_{\mbox{\scriptsize
D}}},\mathbb{R}^3)=\bigl\{v\in
\mathsf{H}^1(\Omega,\mathbb{R}^3)\,\,|\, \,v_{|\Gamma}=0\bigr\}\,,\\
\nonumber \SFQ&=&\mbox{H}_0(\mbox{Curl};\,\Omega,\,\Gamma,\sL(3))\quad\mbox{ defined in (\ref{spacep-bc})}\,,\\
 \nonumber\Lambda&=& L^2(\Omega)\,,
\end{eqnarray} equipped with the norms
\begin{eqnarray}
\nonumber&&  \norm{v}_V:=\norm{\nabla v}_{L^2},\label{norm-V-new}\quad\qquad
\norm{q}_{Q}:=\norm{q}_{\mbox{\scriptsize H}(\mbox{\scriptsize
Curl};\Omega)},\\
&&\norm{z}^2_{Z}:=\norm{v}^2_{V} +\norm{q}^2_{Q}
+\norm{\eta}^2_{L^2}\quad\mbox{ for } z=(v,q,\eta)\in \SFZ\,.
\label{norm-Z-new}
\end{eqnarray}
Now, for the existence of a solution to the problem (\ref{wf-new1})
following Theorem \ref{han-reddy-theo}), we only need
to check that the bilinear form $\ba$ is coercive in $\SFW$.
Following the coercivity inequality obtained  in (\ref{coercive-a}), we immediately get a
positive constant $C= C(m_0,\mu,\alpha_1,L_c,\Omega)>0$ such that
$$a(z,z)\geq C\bigl[\norm{v}_V^2+\norm{\sym q}^2 +\norm{\Curl q}^2_2+\norm{\eta}_2^2\bigr]\,.$$
But this estimate is not enough to establish coercivity. Indeed, the
skew-symmetric (spin) part $\skew q$ of $q$ is not controlled
locally. \vskip.2truecm\noindent
 Motivated by the well-posedness question for
precursors to this model \cite{NCA, EBONEFF}, Neff et al.
\cite{NPW2011-1, NPW2012-1, NPW2012-2, NPW2014}, derived a new
inequality extending Korn's first inequality to incompatible tensor
fields, namely there exists a constant $C(\Omega)>0$ such that
\begin{align}
\label{incompatible_korn}
\forall \, p\in \SFH(\Curl;\,\Omega,\,\BBR^{3\times 3}), \quad & p\times\,n|_\Gamma=0:   \\
  & \underbrace{\|p\|_{L^2(\Omega)}}_{\text{plastic distortion}}\le C(\Omega)\,
     \Big( \underbrace{\norm{\sym p}_{L^2(\Omega)}}_{\text{plastic strain}}+
     \underbrace{ \norm{\Curl p}_{L^2(\Omega)}}_{\text{dislocation density}} \Big)\, .\notag
\end{align}
This shows that if we consider the closure
$\SFH_{\mbox{\scriptsize sym }}(\Curl,\,\Omega,\Gamma;\sL(3))$ of
the linear subspace
$$\{q\in C^\infty(\overline{\Omega},\BBR^{3\times
3})\,\,|\,\,\tr{q}=0,\,\,q\times\,n=0\mbox{ on }\Gamma\}$$ in the
norm
\begin{equation}\label{newnorm-q}\norm{q}_{\mbox{\scriptsize
sym,\,curl}}^2:=\norm{\sym
q}^2_{L^2}+\norm{\mbox{Curl\,}q}^2_{L^2}\,,
\end{equation}
then we have the decisive identity
$$
\SFH_{\mbox{\scriptsize sym
}}(\Curl,\,\Omega,\Gamma;\sL(3))\,\equiv\,
\SFH_0(\Curl,\,\Omega,\Gamma;\sL(3))$$ with equivalence of norms.
Therefore, we have the coercivity inequality

\begin{equation}\label{coercive-a-sym}
\ba(z,z)\geq C\bigl[\norm{v}_V^2+\norm{q}^2_{Q}
+\norm{\eta}_2^2\bigr]=\norm{z}^2_Z\qquad\forall z\in\SFW\,,
\end{equation}
from which we obtain the existence of a solution $(u,p,\gamma_p)\in
\SFW$ to the problem (\ref{wf-new1}). Now setting a posteriori
\begin{equation}\label{omegap}\omega_p(t,x):=\int_0^t\norm{\skew \dot{p}(s,x)}\,ds\,,\end{equation} it
follows that $(u,p,\gamma_p,\omega_p)$ is a solution to the original
problem (\ref{wf-special}).\hfill\qed
\begin{remark}\label{iso-hardspin-alpha2}{\rm
 Notice again that isotropic hardening in the spin-evolution is not necessary for existence of a solution to the problem and it is not
  connected to the uniqueness question either. In fact, arguing as
  in Section \ref{uniqueness} we get the inequality (\ref{unique-ineq}) with
  $\alpha_2=0$, from which and from the flow law in dual
  form on each $\mathcal{S}_k$ of the yield surface, we deduce the
  uniqueness of $u$, $\sigma$, $p$ and $\gamma_p$ while the
  uniqueness of $\omega_p$ follows from (\ref{omegap}) and from the uniqueness of $p$. Therefore, the strong solution is unique also
   in the case where there is no isotropic hardening in the spin-evolution.}\end{remark}

 \subsection{Is it possible to accommodate the special case $\yieldspino=0$ in our model?}\label{hatsigma-equalzero}
 In Gurtin's visco-plastic model \cite{GURT2004} it is possible to consider $\widehat{\sigma}_0=0$. In our setting, this
  case corresponds to the dissipation function

\begin{equation}\label{diss-funct1-new}
\Delta(q,\eta,\beta):=\left\{\begin{array}{ll}\hskip-.2truecm
\yieldzero\,\norm{\sym q} & \mbox{if }\norm{\sym
q}\leq\eta\quad\mbox{and}\quad
\norm{\skew q}\leq\beta\\\\
\hskip-.2truecm\infty & \mbox{otherwise}\,\end{array}\right.
\end{equation}
and the elastic region
\begin{equation}\label{elast-dom}\mathcal{E}:=\bigl\{\Sigma_p=(\Sigma_E,g_1,g_2)\,\,|\,\,\norm{\dev\,\sym\Sigma_E}-\yieldzero+g_1\leq0\mbox{
and }\norm{\skew\Sigma_E}+g_2\leq0\bigr\}\,.\end{equation} The flow
law in dual form is given in Table
\ref{table:flow-flowrule-special1} below.

\begin{table}[h!]\footnotesize\begin{center}
\begin{tabular}{|ll|}\hline &\qquad\\
{\em Free energy:} & $\frac12\langle\C_{\mbox{\scriptsize
iso}}\bvarepsilon_e,\bvarepsilon_e\rangle+\frac12\mu
\,L^2_c\,\norm{\mbox{Curl}\,\PL}^2+\frac12\mu\,\alpha_1\,|\gamma_p|^2+\frac12\mu\,\alpha_2|\omega_p|^2$\\&\\
{\em Elastic region:} &$
\mathcal{E}:=\Bigl\{(\Sigma_E,g_1,g_2)\,\,|\,\,\norm{\dev\,\sym\Sigma_E}-\yieldzero+g_1\leq0\mbox{
and }\norm{\skew\Sigma_E}+g_2\leq0\Bigr\}$\\&\\
{\em Yield surface:} &
$\partial\mathcal{E}=\mathcal{S}_1\cup\mathcal{S}_2 $\\&\\
{\em where} &
$\mathcal{S}_1:=\Bigl\{(\Sigma_E,g_1,g_2)\,\,|\,\,\norm{\dev\,\sym\Sigma_E}-\yieldzero+g_1=0\Bigr\}$\\&\\
&$\mathcal{S}_2:=\Bigl\{(\Sigma_E,g_1,g_2)\,\,|\,\,\norm{\skew\Sigma_E}+g_2=0\Bigr\}$
\\&\\
 {\em Dissipation function:} &$
\Delta(q,\eta,\beta):=\left\{\begin{array}{ll}\hskip-.2truecm
\yieldzero\,\norm{\sym q} & \mbox{if }\norm{\sym
q}\leq\eta\quad\mbox{and}\quad
\norm{\skew q}\leq\beta\\\\
\hskip-.2truecm\infty & \mbox{otherwise}\,\end{array}\right.$
\\&\\
{\em Flow law in dual form:}
&$\left\{\begin{array}{lll}\sym\dot{\PL}&=&\lambda\,
\dsize\frac{\dev\sym\Sigma_E}{\norm{\dev\sym\Sigma_E}}\,,\quad
\skew\dot{p}=0,\quad
\dot{\gamma}_p=\lambda,\quad\dot{\omega}_p=0\quad\mbox{ on
}\mathcal{S}_1
\\\\
\sym\dot{p}&=&0,\quad\skew\dot{\PL}=\lambda\,\dsize\frac{\skew\Sigma_E}{\norm{\skew\Sigma_E}}\,,\quad
\quad\,\,\,\dot{\gamma}_p=0,\quad\dot{\omega}_p=\lambda\quad\mbox{ on
}\mathcal{S}_2
\end{array}\right.$\\&
\\
 \hline
\end{tabular}\caption{\footnotesize The flow rule in dual form in the case
 $\yieldspino=0$ and $\alpha_2>0$.}
   \label{table:flow-flowrule-special1}\end{center}\end{table}
Table \ref{table:flow-flowrule-special1}  illustrates  why both
initial yield stresses $\yieldzero$ and $\yieldspino$ have to be
strictly positive. In fact, since $g_2$ may be zero initially, the
elastic domain $\mathcal{E}$ in (\ref{elast-dom}) may not have
non-empty interior. Therefore, this forbids the use of
$\yieldspino=0$. In the visco-plastic setting
   $\yieldspino=0$ may be accommodated.

\section{The limit case of vanishing characteristic length scale $L_c\to0$}\label{limit-Lc}
 In the limit case $L_c\to0$, looking at the flow rule in its dual formulation, we first observe that the thermodynamic
  driving stress
  $\Sigma_E\in\mathbb{R}^{3\times 3}$ reduces to the symmetric Cauchy stress tensor
  $\sigma\in\mbox{Sym}(3)$ and we see clearly that we do not have the branch $\mathcal{S}_2$ and moreover,

 \begin{equation}\label{limit-flow-S1}
\mbox{on }\mathcal{S}_1:\quad\left\{\begin{array}{ll}\sym\dot{p}\,\,=\lambda\,
\dsize\frac{\dev\sigma}{\norm{\dev\sigma}},&\quad\dot{\gamma}_p=\lambda=\norm{\sym\dot
p},\\ \skew\dot{p}=0,&\quad\dot{\omega}_p=0\,,\end{array}\right.
 \end{equation} while on $\mathcal{S}_3$ we get from the rate-explicit dual formulation

$$ \left\{\hskip-.1truecm\begin{array}{ll}\sym\dot{\PL}=\dsize\frac{2\lambda}{\yieldzero^2}(\norm{\dev\sym\Sigma_E}+g_1)
\dsize\frac{\dev\sym\Sigma_E}{\norm{\dev\sym\Sigma_E}}\,, &\,\,
\dot{\gamma}_p=\dsize\frac{2\lambda}{\yieldzero^2}(\norm{\dev\sym\Sigma_E}+g_1)=\norm{\sym\dot{\PL}}
\\\\
\hskip-.1truecm\skew\dot{\PL}=\dsize\frac{2\lambda}{\yieldspino^2}(\norm{\skew\Sigma_E}+g_2)\dsize\frac{\skew\Sigma_E}{\norm{\skew\Sigma_E}}\,,&\,\,
\dot{\omega}_p=\dsize\frac{2\lambda}{\yieldspino^2}(\norm{\skew\Sigma_E}+g_2)=\norm{\skew\dot{\PL}}
\end{array}\right.$$
in the case $L_c>0$ and $\alpha_2>0$ and

$$\left\{\hskip-.1truecm\begin{array}{ll}\sym\dot{\PL}=\dsize\frac{2\lambda}{\yieldzero^2}(\norm{\dev\,\sym\,\Sigma_E}+g_1)
\dsize\frac{\dev\sym\Sigma_E}{\norm{\dev\sym\Sigma_E}}\,,&\,\,
\dot{\gamma}_p=\dsize\frac{2\lambda}{\yieldzero^2}(\norm{\dev\,\sym\,\Sigma_E}+g_1)=\norm{\sym\dot{\PL}}\\\\
\hskip-.1truecm\skew\dot{\PL}=\dsize\frac{2\lambda}{\yieldspino^2}\skew\Sigma_E\,,&\,\,
\dot{\omega}_p=\dsize\frac{2\lambda}{\yieldspino^2}\norm{\skew\,\Sigma_E}=\norm{\skew\dot{\PL}}
\end{array}\right.$$ in the case $L_c>0$ and $\alpha_2=0$ that altogether
\begin{equation}\label{limit-flow-S3}
\mbox{on }\mathcal{S}_3:\quad\left\{\begin{array}{ll}\sym\dot{p}\,\,=\dsize\frac{2\lambda}{\sigma_0}\,
\dsize\frac{\dev\sigma}{\norm{\dev\sigma}},&\qquad\dot{\gamma}_p=\dsize\frac{2\lambda}{\sigma_0}=\norm{\sym\dot
p},\\\\ \skew\dot{p}=0,&\qquad\dot{\omega}_p=0\,.\end{array}\right.
 \end{equation}
Therefore, we obtain for $L_c\to0$ that all driving stress-tensor quantities are
 symmetric such that, if $p(0)\in\mbox{Sym}(3)$, then we will have
 $p(t)\in\mbox{Sym}(3)$ along the plastic evolution. In that case, our new
 model turns into
 \begin{equation}\label{classical1}\dot{\bvarepsilon}_p=\widehat{\lambda}\,\frac{\dev\sigma}{\norm{\dev\sigma}}\,,\qquad
 \dot{\gamma}_p=\widehat{\lambda}=\norm{\dot{\bvarepsilon}_p}\,,\end{equation}
 which is the dual formulation of the flow rule for classical
 plasticity with isotropic hardening based only on the accumulated
 equivalent plastic strain
 $\gamma_p=\int_0^t\norm{\dot{\bvarepsilon}_p}\,ds$.\\
 For us it is interesting to remark that the evolution of plastic
 spin in our model is related solely to the energetic length scale
 $L_c>0$.
\section{Conclusions and outlook}\label{CO}
From a modelling perspective, it is not difficult to extend the present model to visco-plasticity.
However, the well-posedness result (which we expect to hold) needs to be derived along different methods.
 Moreover, it would be interesting to treat the dynamic case. Both questions are subject of ongoing work.\\
 Since we did not establish unqualified uniqueness in our model (it hinges on the additional regularity
 $\Curl\Curl p\in L^2(\Omega,\mathbb{R}^{3\times 3})$) it will also be interesting to establish higher regularity
  provided the data are regular. It remains open whether we really could have non-uniqueness of the weak solutions if
  regularity is missing. Is the dislocation energy contribution $\Curl p\in L^2(\Omega,\mathbb{R}^{3\times 3})$
   strong enough to prevent non-uniqueness? The question we have to answer is, what least amount of hardening will lead to existence and uniqueness
    in rate-independent gradient plasticity?
\\
  We expect furthermore that a computational implementation suggests itself along the lines of \cite{NSW2009}. Attendant
   to these research perspectives, one should look at simple settings of boundary value problems like anti-plane shear
   to gain more insight in the response of the model and the new features offered by incorporating plastic spin.

  Finally, a major challenge from the mathematical point of view is the replacement of the defect energy  in (\ref{free-eng-expr})$_2$ by a more physically realistic term. The one-homogeneous term 
    $\mu\,L_c\norm{\Curl p}$ was proposed  in \cite{OHNOKUM2007} in the context of single crystal gradient plasticity    and  is
    summarized
     in \cite[p.92]{Han-ReddyBook} while energies of logarithmic form were used in  \cite{FORGUEN2013,Bardepante2015}. However, such defect energies cannot be adopted in the current mathematical framework and hence, their mathematical treatment needs fundamentally new ideas.

Another interesting challenge is the one of studying possible visco-plastic regularizations of our model through either the classical power-law rate-dependence or the proposal in \cite{PANTEBARDE2016}.

\addcontentsline{toc}{section}{Acknowledgements}
\section*{Acknowledgements:}
{ The research of Francois
Ebobisse has been supported by the National Research Foundation
(NRF) of South Africa through the Incentive Grant for Rated
Researchers and the International Centre for Theoretical Physics
(ICTP) through the Associateship Scheme. The first draft of this
work was written at Essen (Germany) while Francois
Ebobisse was visiting the Faculty of Mathematics of the University
of Duisburg-Essen.}

\end{document}